\documentclass[preprint,12pt]{elsarticle}
\usepackage{amsfonts,amssymb,amsmath,mathrsfs}
\usepackage{lineno,hyperref}
\modulolinenumbers[5]

\usepackage{color}

\journal{Arxiv}

\newtheorem{theorem}{Theorem}[section]
\newtheorem{lemma}[theorem]{Lemma}

\newtheorem{proposition}[theorem]{Proposition}
\newtheorem{remark}[theorem]{Remark}

\newtheorem{definition}[theorem]{Definition}
\numberwithin{equation}{section}

\newcommand{\N}{\mathbb{N}}
\newcommand{\Z}{\mathbb{Z}}

\newcommand{\R}{\mathbb{R}}
\newcommand{\C}{\mathbb{C}}
\newcommand{\TT}{\mathbb{T}}
\newcommand{\LL}{\mathbb{L}}

\newenvironment{proof}[1][\noindent \textbf{Proof: }]{#1}{ \hfill $\square$ \vspace{2mm}}

\begin{document}

\begin{frontmatter}
	
	\title{Global Hypoellipticity for Systems in \\ Time-Periodic Gelfand-Shilov Spaces}


	\author[addressUFPR]{Fernando de \'Avila Silva}
	\ead{fernando.avila@ufpr.br}

\author[addressUNITO]{Marco Cappiello}
	\ead{marco.cappiello@unito.it}

\author[addressUFPR]{Alexandre Kirilov}
\ead{akirilov@ufpr.br}

	\address[addressUFPR]{Department of Mathematics, Federal University of Paran\'a, Caixa Postal 19081, \\ CEP 81531-980, Curitiba, Brazil}

	\address[addressUNITO]{Department of Mathematics, University of Turin, Via Carlo Alberto 10, 10123, Turin, Italy}

	\begin{abstract}
	We investigate the global hypoellipticity of a class of overdetermined systems with coefficients depending both on time and space variables in the setting of time-periodic Gelfand-Shilov spaces. Our main result provides necessary and sufficient conditions for the global hypoellipticity of this class of systems, stated in terms of Diophantine-type estimates and sign-changing behavior of the imaginary parts of the coefficients. Through a reduction to a normal form and detailed construction of singular solutions, we fully characterize when the system fails to be globally hypoelliptic. 	
\end{abstract}
	

\end{frontmatter}

\section{Introduction and Main Result}

The aim of this paper is to investigate the global hypoellipticity for a class of overdetermined systems of evolution equations. This mathematical problem has a long hystory. First of all, we recall that the \emph{local} hypoellipticity and solvability of systems are described  in the seminal work of Treves \cite{Tre76}. Also the global investigations for differential operators on the compact setting, with special attention to vector fields on the torus, count many contributions. In this context, we highlight the paper \cite{BerCorMal93} by Bergamasco, Cordaro and Malagutti. We point out that these investigations are also developed in the setting of analytic and Gevrey spaces (see for instance  \cite{AriKirMed19,Ber99} and the references quoted therein). 
\\ \indent		
In this paper we consider systems of the type 
\begin{equation}\label{general-const-system}
	L_r u = f_r, \qquad r = 1, \ldots, m,
\end{equation}
where
\begin{equation}\label{general-const-system2}    
	L_r = D_{t_r} + (a_r + ib_r)(t_r)P(x, D_x), \qquad r = 1, \ldots, m,
\end{equation}
and the real-valued coefficients $a_r$ and $b_r$ belong to the Gevrey classes (of order $\sigma > 1$) on the torus. The functions $f_r$ and the distribution $u$ are taken in the Gelfand-Shilov classes on $\TT^m \times \R^n$, as introduced in \cite{AviCap22}. The operator $P(x, D_x)$ is a normal differential operator of the form
\begin{equation}\label{P-intro1}
	P = P(x, D_x) = \sum_{|\alpha| + |\beta| \leq M} c_{\alpha, \beta} x^{\beta} D_x^{\alpha}, \quad c_{\alpha, \beta} \in \C,
\end{equation}
of order $M \geq 2$, which satisfies the global ellipticity condition 
\begin{equation}\label{P-elliptic-intro}
	p_M(x, \xi) = \sum_{|\alpha| + |\beta| = M} c_{\alpha, \beta} x^{\beta} \xi^{\alpha} \neq 0, \quad (x, \xi) \neq (0, 0).
\end{equation}

The normality condition and \eqref{P-elliptic-intro} imply that $P$ has a discrete spectrum consisting of a sequence of real eigenvalues $\lambda_j$, with $|\lambda_j| \to +\infty$ as $j \to +\infty$. The corresponding eigenfunctions $\varphi_j$ form an orthonormal basis in $L^2(\R^n)$. Furthermore, the following asymptotic Weyl formula holds:
\begin{equation*}
	|\lambda_j| \sim \rho  j^{\, M / 2n}, \quad \text{as } j \to \infty,
\end{equation*}
for some positive constant $\rho$.
\\ \indent Elliptic equations involving operators of the form \eqref{P-intro1}, \eqref{P-elliptic-intro} and their pseudodifferential extensions have been studied in \cite{CapGraRod06, CapGraRod10jam, Shu87}, both in the Schwartz setting and and within the framework of Gelfand-Shilov spaces $\mathcal{S}^{\mu}_\nu(\R^n)$ defined as the spaces of all smooth functions $f$ on $\R^n$ satisfying the following condition:
\begin{equation}
	\sup_{x \in \R^n}\sup_{\alpha \in \N^n} C^{-|\alpha|} (\alpha!)^{-\mu} \exp \left(c|x|^{1/\nu} \right) |\partial^\alpha f(x)| <\infty
\end{equation}
for some constants $C,c>0$. These spaces represent a natural global counterpart on $\R^n$ of the Gevrey classes and have acquired increasing relevance over the last two decades in the analysis of partial differential equations with Gevrey regular coefficients, see \cite{Ari24,AriAscCap22,AscCap08,AscCap19,CapGraRod06,CapGraRod10jam,CapGraRod10cpde,CapPilPra16,CorNicRod15rmi,CorNicRod15tams,NicRod10book}. Indeed, Gelfand-Shilov spaces are the most natural setting for the analysis of the operators \eqref{P-intro1}, \eqref{P-elliptic-intro} since the eigenfunctions of these operators are indeed in $\mathcal{S}^{1/2}_{1/2}(\R^n)$ and the same holds for the solutions of the equation $Pu=0$.
More recently in \cite{AviCap22, AviCap24}, the first two authors obtained necessary and sufficient conditions for global hypoellipticity and for global solvability for evolution operators of the form
\begin{equation}\label{scalar-operator}    
	L = D_{t} + (a+ib)(t)P(x, D_x), \qquad (t,x) \in \TT^1 \times \R^n, 
\end{equation}  
in the frame of time-periodic Gelfand-Shilov spaces $\mathcal{S}_{\sigma, \mu}(\TT^1 \times \R^n)$, with $\sigma \geq 1$ and $\mu \geq 1/2$. For $m \in \N$, the space $\mathcal{S}_{\sigma, \mu}(\TT^m \times \R^n)$ (denoted simply as $\mathcal{S}_{\sigma, \mu}$), is defined as the set of all smooth functions $u$ on $\TT^m \times \R^n$ satisfying
\begin{equation}\label{firstnorm}
	|u|_{\sigma, \mu, C} :=
	\sup_{\alpha, \beta \in \N^n, \gamma \in \N^m} C^{-|\alpha + \beta| - |\gamma|} \gamma!^{-\sigma} (\alpha! \beta!)^{-\mu} \sup_{(t, x) \in \mathbb{T}^m \times \R^n} |x^\alpha \partial_x^\beta \partial_t^\gamma u(t, x)|
\end{equation}
is finite for some positive constant $C$. The same has been done in \cite{AviCapKiri24} for systems of the form
\begin{equation}\label{general-const-system-paper1}
	L_r u=f_r, \qquad r=1,\ldots,m,
\end{equation} with
\begin{equation}\label{general-const-system2-paper1}	L_r = Q_r(D_t)+ d_rP(x,D_x), \qquad r=1,\ldots,m,
\end{equation}
where $d_r \in \C$, and $f_r$ are assigned functions or distributions and the operator
$$
Q_r(D_t)=\sum_{|\alpha|\leq k_r} c_{\alpha,r} D_t^\alpha
$$
is a constant-coefficient differential operator of order $k_r$ on the torus $\TT^m$, with $\alpha= (\alpha_1, \ldots, \alpha_n) \in \mathbb{N}_0^m$,  $c_{\alpha,r} \in \C,$ and $D_t^\alpha = D_1^{\alpha_1} D_2^{\alpha_2}\cdots D_m^{\alpha_m}$.

In the present paper we aim to extend our previous results and consider systems with time depending coefficients of the form \eqref{general-const-system}, \eqref{general-const-system2}.

Throughout this article, we consider the space
$$
\mathscr{F}_\mu := \bigcup_{\sigma > 1} \mathcal{S}_{\sigma, \mu}
$$
endowed with the topology
$$
\mathscr{F}_\mu = \underset{\sigma \to +\infty}{\operatorname*{ind\, lim}} \;\mathcal{S}_{\sigma, \mu}.
$$
Thus, $\mathscr{F}_\mu$ is an inductive limit of Banach spaces. We denote by $\mathscr{F}'_\mu$ the dual space of continuous linear functionals on $\mathscr{F}_\mu$. Further details on the topology of $\mathscr{F}_\mu$ and $\mathscr{F}'_\mu$ will be provided in Subsection \ref{FandF'}.
Next, we give a suitable definition of global hypoellipticity for systems of the type \eqref{general-const-system}, \eqref{general-const-system2}.

\begin{definition}\label{defGH}
	System \eqref{general-const-system} is said to be $\mathscr{F}_\mu$-globally hypoelliptic if, whenever $u \in \mathscr{F}'_\mu$ and $L_r u \in \mathscr{F}_\mu$ for $1 \leq r \leq m$, it follows that $u \in \mathscr{F}_\mu$.
\end{definition}

\begin{remark} \label{GHcomparison} Note that the notion of global hypoellipticity given in the above definition is more restrictive than the one used in \cite{AviCap22} for the case $m = 1$. In fact, in \cite{AviCap22} $u$ is assumed a priori to be in the set $\mathscr{U}_\mu = \bigcup_{\sigma > 1} \mathcal{S}_{\sigma, \mu}'$, which contains $\mathscr{F}'_\mu$ and the operator $$\mathcal{L}=D_t+(a(t)+ib(t))P$$  is said to be globally hypoelliptic (GH) if $\mathcal{L} u \in \mathscr{F}_{\mu}(\mathbb{T} \times \R^n)$ and $u \in \mathscr{U}_{\mu}(\mathbb{T} \times \R^n)$ imply $u \in \mathscr{F}_{\mu}(\mathbb{T} \times \R^n)$. 
	Now, if $\mathcal{L}$ is GH and 
	$$
	u \in \mathscr{F}'_\mu(\mathbb{T} \times \R^n) \subset \mathscr{U}_\mu(\mathbb{T} \times \R^n)
	$$
	is a solution of $\mathcal{L} u \in \mathscr{F}_{\mu}(\mathbb{T} \times \R^n)$, then we get $u \in \mathscr{F}_\mu(\mathbb{T} \times \R^n)$. Therefore, $\mathcal{L}$ is $\mathscr{F}_\mu$-globally hypoelliptic according to Definition \ref{defGH}. This shows that
	\begin{equation}\label{GHimpliesFGH}
		\text{global hypoellipticity } \Longrightarrow \  \mathscr{F}_\mu\text{-global hypoellipticity}.
	\end{equation}
	The modification of the notion of global hypoellipticity adopted here is necessary because the arguments in this paper rely on the reduction of the system $\mathbb{L}$ to its normal form via a conjugation, which possesses convenient mapping properties on $\mathscr{F}'_\mu$, but not on $\mathscr{U}_\mu$. Nevertheless, for the sake of consistency, in the sequel we shall reformulate the results obtained in \cite{AviCap22} in terms of the notion of $\mathscr{F}_\mu$-global hypoellipticity, see Proposition \ref{GHm=1}.
\end{remark}

Before stating our main theorem, we introduce a notion of vector approximation in $\mathbb{R}^n$ that resembles the concept of Liouville exponential vectors of order $\sigma$, which typically appear in the study of global hypoellipticity for differential operators on the $m$-dimensional torus, see \cite{Ber99,Pet05}.

\begin{definition}
	Let $\sigma>1$ and $\mu\geq 1/2$ and let $\{\lambda_\ell\}_{\ell \in \N}$ be a sequence of real numbers. We say that a vector  $\alpha=(\alpha_1,\ldots,\alpha_N)$ in $\mathbb{R}^N$ satisfies condition $\mathscr{D}_{\sigma,\mu}$ with respect to the sequence $\{\lambda_\ell\}_{\ell \in \N}$ if the following property holds true: for every $\epsilon>0$ there exists $C_{\epsilon}>0$ such that
	\begin{equation}\label{non-exp-Liou}
		\| \tau - \alpha \lambda_\ell\| \geqslant C_{\epsilon} \exp \left(-\epsilon(\|\tau\|^{1/\sigma}  + \ell^{1/(2n\mu)})\right),
	\end{equation}
	for all $(\tau,\ell) \in \Z^N \times \N$ such that $\tau - \alpha \lambda_\ell \neq 0$.
\end{definition}

Now, for each  $r\in \{1, \dots, m\}$, we define
\begin{equation}\label{realpartaverage}
	a_{r,0} := (2 \pi)^{-1} \int_{0}^{2\pi} a_{r}(s)\, ds. 
\end{equation}

\begin{theorem}\label{main_Theorem}
	Let $\mu \geq 1/2.$ The	system \eqref{general-const-system}  is $\mathscr{F}_{\mu}$-globally hypoelliptic if and only if at least one of the following conditions occurs:	
	\begin{enumerate}
		\item[(I)] There is $r\in \{1, \dots, m\}$ such that the function $\mathbb{T} \ni t_r  \mapsto b_r(t_r)$ does not change sign and is not identically zero;
		\item[(II)] If the set
		\begin{equation}\label{J}
			J = \left\{ r\in \{1, \dots, m\}: b_r(t_r) \equiv 0 \right\} ,
		\end{equation}
		is not empty and $J =\{r_1< \dots< r_k \}$,
		then $a_{J0} = (a_{r_1, 0}, \dots, a_{r_k, 0})$ satisfies condition $\mathscr{D}_{\sigma,\mu}$ with respect to the sequence $\{\lambda_\ell\}_{\ell \in \N}$ of the eigenvalues of $P$ for all $\sigma\geq M\mu$.
	\end{enumerate}
\end{theorem}

Let us present an example of a vector satisfying condition $\mathscr{D}_{\sigma,\mu}$ for all $\sigma > 1$. Consider first the Schrödinger  harmonic oscillator operator
\[
H = -\Delta + \|x\|^2,
\]
defined on $\mathbb{R}^n$, whose eigenvalues are given by
\[
\lambda_k = \sum_{j=1}^{n} (2k_j + 1), \quad k = (k_1, \ldots, k_n) \in \mathbb{N}_0^n.
\]

We recall that a vector $\alpha \in \mathbb{R}^n \setminus \mathbb{Q}^n$ is said to be exponentially Liouville of order $\sigma \geq 1$ (simply exponentially Liouville when $\sigma = 1$) if there exists $\epsilon > 0$ such that the inequality
\[
\|\tau - \ell \alpha\| < \exp\left(-\epsilon |\ell|^{1/\sigma}\right)
\]
admits infinitely many solutions $(\tau, \ell) \in \mathbb{Z}^n \times \mathbb{Z}$. Consequently, $\alpha$ is not exponentially Liouville of order $\sigma \geq 1$ if, for every $\epsilon > 0$, there exists a constant $C_\epsilon > 0$ such that
\begin{equation}\label{Berg-vector-nL}
	\|\tau - \ell \alpha\| \geq C_\epsilon \exp\left(-\epsilon |\ell|^{1/\sigma}\right),
\end{equation}
for all \((\tau, \ell) \in \mathbb{Z}^n \times \mathbb{Z}\), satisfying \( \tau - \ell \alpha \neq 0.\)

We claim that if $\alpha$ is not exponentially Liouville of order $\sigma$ for every $\sigma \geq 2\mu$ with respect to the eigenvalues $\lambda_k$, then it satisfies $\mathscr{D}_{\sigma,\mu}$ for all $\sigma \geq 2\mu$. Indeed, if this were not the case, there would exist some $\sigma \geq 2\mu $, $\epsilon > 0$, and a sequence $(\tau_k, j_k) \in \mathbb{Z}^n \times \mathbb{N}$ such that
\[
0 < \|\tau_k - \alpha \lambda_{j_k}\| < \exp\left(-\epsilon(\|\tau_k\|^{1/\sigma} + j_k^{1/(2n\mu)})\right) \leq \exp\left(-\epsilon j_k^{1/(2n\mu)}\right).
\]
Since $\lambda_{j_k} \sim c\, j_k^{1/n}$ for some constant $c > 0$, it follows that
\[
0 < \|\tau_k - \alpha \lambda_{j_k}\| < \exp\left(-\epsilon' \lambda_{j_k}^{1/(2\mu)}\right)\leq \exp\left(-\epsilon' \lambda_{j_k}^{1/\sigma}\right),
\]
where $\epsilon' = \epsilon c^{-1/(2\mu)}$, which contradicts \eqref{Berg-vector-nL} for $\sigma \geq 2\mu$.

\medskip

The paper is structured as follows. In Section \ref{sec-notations} we describe the time-periodic Gelfand-Shilov spaces introduced in \cite{AviCap22} and present a characterization  based on their Fourier coefficients. Moreover, in Section \ref{FandF'} we introduce the spaces  $\mathscr{F}_\mu$ and $\mathscr{F}'_\mu$ and their topologies. 
In Section \ref{secRedNormalForm} we exhibit the normal form of system  \eqref{general-const-system}, namely, the suitable system \eqref{L_0} where the real parts of the coefficients are time-independent  (see Theorem \ref{The-Normal}). In particular, it follows from Proposition \ref{Prop_real_valued_coeff} that the systems involving only real-valued coefficients can be reduced to the time-independent case. Finally, the proof of Theorem \ref{main_Theorem} is developed in Sections \ref{secSuffCOnditions} (sufficiency) and \ref{necessity-section} (necessity).

\section{Fourier Analysis in Gelfand-Shilov spaces  \label{sec-notations}} 

In this section, we introduce the essential tools of Fourier analysis necessary for a precise presentation of our results and examples, as well as their proofs. Although some of the techniques discussed here have been explored in greater detail in earlier works \cite{AviCap22,AviCap24,AviCapKiri24}, we have developed new methods specifically for this study, that will also be presented.

\medskip
Let $\mathcal{G}^{\sigma,h}(\TT^m)$, with $h > 0$ and $\sigma \geq 1$, denote the space of all smooth functions $\varphi \in C^{\infty}(\TT^m)$ for which there exists a constant $C > 0$ such that
\begin{equation*}
	\sup_{t \in \TT^m} |\partial^{\gamma}\varphi(t)| \leq C h^{\gamma}(\gamma!)^{\sigma}, \quad \gamma \in \N_0^m.
\end{equation*}

This space is a Banach space with respect to the norm
\begin{equation*}
	\|\varphi\|_{\sigma,h} = \sup_{\gamma \in \N^m_0}\left\{\sup_{t \in \TT^m} |\partial^{\gamma}\varphi(t)| h^{-\gamma}(\gamma!)^{-\sigma}\right\}.
\end{equation*}

The space of periodic Gevrey functions of order $\sigma$ is defined as
\[
\mathcal{G}^{\sigma}(\TT^m) = \underset{h\rightarrow +\infty}{\operatorname*{ind\,lim\,} } \mathcal{G}^{\sigma, h} (\TT^m),
\]
and its dual space is denoted by $(\mathcal{G}^{\sigma})'(\TT^m)$.

Gevrey functions $f \in \mathcal{G}^{\sigma}(\TT^m)$ and ultradistributions $u \in (\mathcal{G}^{\sigma})' (\TT^m)$ can be represented by their Fourier series:
\begin{equation*}
	f(t) = \sum_{\tau \in \Z^m} \widehat{f}(\tau) e^{i \tau \cdot t}, \quad \text{and} \quad u(t) = \sum_{\tau \in \Z^m} \widehat{u}(\tau) e^{i \tau \cdot t},
\end{equation*}
where the Fourier coefficients are given by
\begin{equation*}
	\widehat{f}(\tau) = \frac{1}{(2\pi)^m} \int_{\TT^m} f(t) e^{-i \tau \cdot t} \, dt, \quad \text{and} \quad \widehat{u}(\tau) = \frac{1}{(2\pi)^m} \langle u, e^{-i \tau \cdot t} \rangle,  \quad \tau \in \Z^m.
\end{equation*}

The decay or growth rate of the Fourier coefficients is used to characterize functions and distributions in the following way:
\begin{itemize}
	\item[\(\circ\)] $f \in \mathcal{G}^{\sigma}(\TT^m)$ if and only if there exist constants $\varepsilon, C > 0$ such that
	\begin{equation}\label{seq-coeff-gevrey-func}
		|\widehat{f}(\tau)| \leq C e^{-\varepsilon\|\tau\|^{\frac{1}{\sigma}}},  \quad \forall \tau \in \Z^m;
	\end{equation}
	\item[\(\circ\)] $u \in (\mathcal{G}^{\sigma})'(\TT^m)$ if and only if for all $\varepsilon > 0$, there is a constant $C_\varepsilon > 0$ such that
	\begin{equation}\label{seq-coeff-gevrey-distrib}
		|\widehat{u}(\tau)| \leq C_\varepsilon e^{\varepsilon\|\tau\|^{\frac{1}{\sigma}}}, \quad \forall \tau \in \Z^m.
	\end{equation}
\end{itemize}

\subsection{The spaces $\mathcal{S}_{\sigma, \mu}$ and $\mathcal{S}_{\sigma, \mu}'$}

The space $\mathcal{S}_{\sigma, \mu} (\mathbb{T}^m \times \R^n)$ was introduced in Section 1 (see \eqref{firstnorm}), and $\mathcal{S}'_{\sigma, \mu} (\mathbb{T}^m \times \R^n)$ denotes its dual. Here, we revisit key results related to their Fourier analysis, as detailed in \cite{AviCap22} and \cite{AviCapKiri24}. In addition, we introduce the partial Fourier series with respect to selected time variables, which play a crucial role in developing our main results. This approach allows us to fully characterize functions and ultradistributions in terms of their partial Fourier series.

Let $\varphi_j \in \mathcal{S}_{^{1/2}}^{_{1/2}} (\R^n)$ be the eigenfunctions of the operator $P$ defined in \eqref{P-intro1}, with $j \in \N$. If $u \in \mathcal{S}_{\sigma, \mu}' (\mathbb{T}^m \times \R^n)$, then for any $j \in \N$, the linear form $u_j: \mathcal{G}^{\sigma}(\TT^m) \to \C$ given by
\begin{equation}\label{partial} 
	\langle u_j(t), \psi(t) \rangle \doteq \langle u(t, x), \psi(t)\varphi_j(x) \rangle 
\end{equation}
belongs to $(\mathcal{G}^{\sigma})'(\TT^m)$. Moreover, for each $\varepsilon > 0$ and $h > 0$, there exists a constant $C_{\varepsilon, h} > 0$ such that for all $j \in \N$ and $\psi \in \mathcal{G}^{\sigma, h}(\TT^m)$, we have
\begin{equation} 
	|\langle u_j(t), \psi(t) \rangle| \leq C_{\varepsilon, h} \| \psi \|_{\sigma, h} \exp \left(\varepsilon j^{\frac{1}{2n\mu}}\right).
\end{equation}

Hence, the action of $u$ on $\mathcal{S}_{\sigma, \mu}$ is given by
\begin{equation}\label{def_serie}
	\langle u, \theta \rangle  =
	\sum_{j \in \N} \langle u_j(t)\varphi_j(x), \theta(t, x) \rangle,
\end{equation}
where
\begin{equation}\label{uj-phij-theta-action}
	\langle u_j(t)\varphi_j(x), \theta \rangle \doteq 
	\left\langle u_j(t), \int_{\mathbb{R}^n} \theta(t, x) \varphi_j(x) \, dx \right\rangle,
\end{equation}
and $\{u_j(t)\}_{j \in \N}$ is a sequence in $(\mathcal{G}^{\sigma})'(\TT^m)$ as defined in \eqref{partial}.

\medskip
Conversely, given a sequence $\{u_j(t)\}_{j \in \N}$ in $(\mathcal{G}^{\sigma})'(\TT^m)$ such that for every $\varepsilon > 0$ and $h > 0$, there exists a constant $C_{\varepsilon, h} > 0$ satisfying
\begin{equation}\label{estimate-distr}
	|\langle u_j, \psi \rangle| \leq C_{\varepsilon, h} \|\psi\|_{\sigma, h} \exp\left(\varepsilon j^{\frac{1}{2n\mu}}\right),
\end{equation}
for all $j \in \N$ and $\psi \in \mathcal{G}^{\sigma, h}(\TT^m)$, we have that 
$$
u = \sum_{j \in \N} u_j(t) \varphi_j(x) \in \mathcal{S}_{\sigma, \mu}'(\mathbb{T}^m \times \R^n),
$$
and
$$
\langle u_j, \psi(t) \rangle = \langle u, \psi(t)\varphi_j(x) \rangle,
$$
for every $\psi \in \mathcal{G}^{\sigma}(\TT^m)$.

\begin{remark}
	According to \cite[Lemma 3.1]{GraPilRod11}, an equivalent norm to \eqref{firstnorm} on $\mathcal{S}_{\sigma, \mu, C}$ is given by:
	\begin{equation}\label{secondnorm}
		\| u \|_{\sigma, \mu, C} := \sup_{\gamma \in \N^m, M \in \N} C^{-M - |\gamma|} M!^{-m\mu} \gamma!^{-\sigma} \| P^M \partial_t^\gamma u \|_{L^2(\TT^m \times \R^n)}.
	\end{equation}
	To take advantage of the good properties of the operator $P(x, D_x)$, in some cases it is convenient to use the norm \eqref{secondnorm} instead of \eqref{firstnorm}, cf. for instance the proof of Proposition \ref{Theorem-Psi} below.
\end{remark}

\begin{theorem}\label{fourier_char_gelfand}
	Let $\mu\geq 1/2$ and $\sigma \geq 1$, and let $u \in \mathcal{S}_{\sigma,\mu}'$. Then $u \in \mathcal{S}_{\sigma,\mu}$ if and only if it can be represented as 
	\begin{equation*}
		u(t,x) = \sum_{j \in \N} u_j(t) \varphi_j(x),
	\end{equation*}
	where
	$$
	u_j(t) = \int_{\R^n} u(t,x)\varphi_j(x)dx,
	$$
	and there exist $C >0$ and $\varepsilon>0$ such that
	\begin{equation} \label{deccoeff}
		\sup_{t \in \TT^m} | \partial_t^\gamma u_j(t)| \leq
		C^{|\gamma|+1} (\gamma!)^{\sigma} \exp \left[-\varepsilon j^{\frac{1}{2n\mu}} \right],
	\end{equation}		
	for all $j \in \N$ and  $\gamma \in \N^m$.
\end{theorem}
\begin{proof}
	See	\cite[Theorem 2.4]{AviCap22}.
\end{proof}

Given that the elements in $\mathcal{S}_{\sigma,\mu}$ and $\mathcal{S}_{\sigma,\mu}'$ have Fourier coefficients that lie on Gevrey classes in the torus, we may analyze their Fourier representations. Consequently, we derive Fourier coefficients with respect to all variables $t$ and $x$ simultaneously, referred to here as the full Fourier coefficients of the elements of $\mathcal{S}_{\sigma,\mu}$ and $\mathcal{S}_{\sigma,\mu}'$.

\begin{theorem}\label{charac_full_fourier-functions}
	Let $\{a(\tau,j)\}_{(\tau,j) \in \Z^m \times \N}$ be a sequence of complex numbers, and consider the formal series
	$$
	a(t,x) = \sum_{j \in \N} \sum_{\tau \in \Z^m} a(\tau,j) e^{it\cdot \tau}\varphi_j(x),
	$$
	where $(t,x)\in\TT^m\times\R^n$. Then, 
	\begin{enumerate}
		\item 	$a \in \mathcal{S}_{\sigma,\mu}$ if and only if there exist $\varepsilon, C>0$ such that, for every $(\tau,j) \in \Z^m\times\N$,
		\begin{equation}\label{seq-full-coeff-funct}
			|a(\tau,j)| \leq C 
			\exp\left[-\varepsilon \left(\|\tau\|^{\frac{1}{\sigma}} + j^{\frac{1}{2n\mu}}\right)\right].
		\end{equation}
		
		Moreover, under this condition,  we have $\widehat{a_j}(\tau) = a(\tau,j)$ for all $(\tau,j) \in \Z^m\times \N$, where
		$$
		\widehat{a_j}(\tau)=\frac1{(2\pi)^m}\int_{\TT^m}a_j(t)e^{-it\cdot\tau}dt,  \quad \text{and} \quad  a_j(t) = \int_{\R^n} a(t,x)\varphi_j(x)dx.
		$$
		\item $a \in \mathcal{S}_{\sigma,\mu}'$ if and only if for every $\varepsilon>0$, there exists $C_{\varepsilon}>0$ such that, for every $(\tau,j) \in \Z^m\times\N$,
		\begin{equation}\label{seq-full-coeff-distr}
			|a(\tau,j)| \leq C_{\varepsilon} 
			\exp\left[\varepsilon \left(\|\tau\|^{\frac{1}{\sigma}} + j^{\frac{1}{2n\mu}}\right)\right].
		\end{equation}
		
		Moreover, under this condition,  we have $\widehat{a_j}(\tau) = a(\tau,j)$ for all $(\tau,j) \in \Z^m\times \N$, where
		$$
		\langle  a_j(t) \, , \, \psi(t) \rangle \doteq 
		\langle  a ,  \psi(t)\varphi_j(x) \rangle, \quad \psi \in \mathcal{G}^{\sigma}(\TT^m).
		$$
	\end{enumerate}
\end{theorem}
\begin{proof}
	See	\cite[Theorems 2.4 and 2.5]{AviCapKiri24}.
\end{proof}

Next, to introduce the partial Fourier series with respect to specific coordinates,  let $m = p+q$ and write $\TT^m = \TT^p \times \TT^q$. In particular, we write
$$
t = (t',t'')\in \TT^m , \text{ meaning that } \ t' \in \TT^p  \ \textrm{
	and } \ t'' \in \TT^q.
$$
Given $f \in \mathcal{G}^{\sigma}(\TT^m)$ we can write
$$
f(t',t'') = \sum_{\eta \in \Z^q} \widehat{f}(t',\eta) e^{i t'' \cdot \eta},
$$
where
$$
\widehat{f}(t',\eta) = \dfrac{1}{(2 \pi)^q} \int_{\TT^{q}}
f(t',t'') e^{-i t'' \cdot \eta} dt''  \in \mathcal{G}^{\sigma}(\TT^p), \ \eta \in \Z^q.
$$
We call $\widehat{f}(t',\eta)$ the partial Fourier coefficient of $f$ with respect to the variables $t''$. Thus, any $f \in \mathcal{S}_{\sigma,\mu} (\mathbb{T}^m\times \R^n)$ can be expressed as:
$$
f(t,x) = f(t',t'',x) = \sum_{j \in \N} \sum_{\eta \in \Z^q} \widehat{f}_j(t',\eta) e^{it''\cdot \eta}\varphi_j(x).
$$

Similarly, the definition for distributions is provided in a manner consistent with natural duality principles.

\begin{theorem}\label{charac_partial_fourier-functions}
	Let $\{a_{j,\eta}(t')\}_{(j,\eta) \in \N \times \Z^q}$ be a sequence in $\mathcal{G}^{\sigma}(\TT^p)$, and consider the formal series
	\begin{equation}\label{sumpartialFour} 
		a(t',t'',x) = \sum_{j \in \N} \sum_{\eta \in \Z^q} a_{j,\eta}(t') e^{it''\cdot \eta}\varphi_j(x).
	\end{equation}
	Then, $a \in \mathcal{S}_{\sigma,\mu}$ if and only if there exist
	$\varepsilon>0$ and $C>0$ such that
	\begin{equation}\label{seq-partial-coeff-funct}
		|\partial_{t'}^{\alpha}a_{j,\eta}(t')| \leq C^{|\alpha|} (\alpha !)^{\sigma}
		\exp\left[-\varepsilon \left(\|\eta\|^{\frac{1}{\sigma}} + j^{\frac{1}{2n\mu}}\right)\right], \ t' \in \TT^p,
	\end{equation}
	for all $\alpha \in \N_0^p,  \eta \in \Z^q$ and  $j \in \N$.
	
	Moreover, when this condition holds, we have $\widehat{a_j}(t', \eta) = a_{j,\eta}(t')$ for all $(\tau,j) \in \Z^m\times \N$, where
	$$
	a_j(t) = \int_{\R^n} a(t,x)\varphi_j(x)dx.
	$$
\end{theorem}
\begin{proof}
	Assume that \eqref{seq-partial-coeff-funct} holds. Then we have that
	\begin{eqnarray*}|\partial_{t'}^\alpha \partial_{t''}^\beta a_{j,\eta}(t')e^{it''\eta}| &\leq& C^{|\alpha|+1}\alpha!^\sigma \|\eta\|^{|\beta|} \exp \left[-\varepsilon \left(\|\eta\|^{\frac{1}{\sigma}} + j^{\frac{1}{2n\mu}}\right)\right] \\ &\leq&
		C_\varepsilon^{|\alpha|+|\beta|+1}(\alpha!\beta!)^\sigma  \exp \left[-\frac{\varepsilon}2 \left(\|\eta\|^{\frac{1}{\sigma}} + j^{\frac{1}{2n\mu}}\right)\right]
	\end{eqnarray*}
	
	This is sufficient to conclude that the sum in \eqref{sumpartialFour} converges in $\mathcal{S}_{\sigma, \mu}$. Conversely, if the sum converges in $\mathcal{S}_{\sigma, \mu}$, then there exists $\varepsilon>0$ such that
	$$|\partial_{t'}^\alpha \hat{a}_j(t',\eta)| \leq \frac{C^{|\alpha|+|\gamma|+1}(\alpha!\gamma!)^\sigma}{|\eta|^{|\gamma|}} \exp(-\varepsilon j^{\frac1{2n\mu}}).$$
	Taking the infimum over $\gamma \in \Z^q$, we obtain that
	$$|\partial_{t'}^\alpha \widehat{a}_j(t',\eta)| \leq C_1^{|\alpha|}\alpha!^\sigma \exp(-\varepsilon (\|\eta\|^{\frac1{\sigma}}+ j^{\frac1{2n\mu}})$$
	from which the estimate \eqref{seq-partial-coeff-funct} follows.
\end{proof}

\subsection{The spaces $\mathscr{F}_\mu$ and $\mathscr{F}'_\mu$}\label{FandF'}

Recall that the Gelfand-Shilov classes, for $\sigma > 1$ and $\mu \geq 1/2$, are constructed from the Banach spaces $\mathcal{S}_{\sigma, \mu, C}$, with $C > 0$, consisting of smooth functions $f: \TT^m \times \R^n \to \C$ for which the norm \eqref{firstnorm} is finite. We then equip $\mathcal{S}_{\sigma, \mu} = \bigcup_{C > 0} \mathcal{S}_{\sigma, \mu, C}$ with the inductive limit topology
$$
\mathcal{S}_{\sigma, \mu} = \varinjlim_{C \rightarrow +\infty} \mathcal{S}_{\sigma, \mu, C}.
$$

Furthermore, $\mathcal{S}_{\sigma, \mu}'$ denotes the space of all linear continuous forms $u: \mathcal{S}_{\sigma, \mu} \to \C$. Here, the continuity means that a sequence $\varphi_j \to 0$ in $\mathscr{F}_\mu$ if and only if there exists $\sigma >1$ such that $\varphi_j \in \mathcal{S}_{\sigma, \mu}$ for every $j$ and $\varphi_j \to 0$ in $\mathcal{S}_{\sigma, \mu}$.  This is in turn equivalent to ask if there exists $C>0$ such that $\varphi_j \in \mathcal{S}_{\sigma, \mu,C}$ for every $j$ and that $\varphi_j \to 0$ in $\mathcal{S}_{\sigma, \mu,C}$. Hence a sequence $\varphi_j \to 0$ in $\mathscr{F}_\mu$ if and only if there exist $C>0$ and $\sigma >1$ such that 
$\varphi_j \in \mathcal{S}_{\sigma, \mu,C}$ for every $j$ and $\varphi_j \to 0$ in $\mathcal{S}_{\sigma, \mu,C}$. Equivalently we can write 
\begin{equation}\label{identity}
	\mathscr{F}_\mu = \bigcup_{\sigma >1} \mathcal{S}_{\sigma, \mu, \sigma -1}.
\end{equation}

This allows us to express $\mathscr{F}_\mu$ as an inductive limit of the Banach spaces $\mathcal{S}_{\sigma, \mu,\sigma-1}$. It should be noted that $\sigma -1$ accepts all positive real values since $\sigma >1$. Moreover, it is not difficult to prove the identity \eqref{identity} and the equivalency of the two topologies defined on $\mathscr{F}_\mu$. 

Next, the dual space $\mathscr{F}'_\mu$ of all continuous and linear functionals on $\mathscr{F}_\mu$ may be described in the following way.

\begin{proposition}
	A linear form $u : \mathscr{F}_\mu \to \C$ is an element of $\mathscr{F}'_\mu$ if and only if for every $C > 0$ and $\sigma > 1$, there exists a constant $B_{C, \sigma} > 0$ such that
	\begin{equation}\label{equivdual}
		|\langle u , \varphi\rangle| \leq B_{C, \sigma} \sup_{\alpha, \beta, \gamma} \ \sup_{(t, x) \in \TT^m \times \R^n} C^{-|\alpha| - |\beta| - |\gamma|} (\alpha! \beta!)^{-\mu} (\gamma!)^{-\sigma} |x^\alpha \partial_x^\beta \partial_t^\gamma \varphi(t, x)|
	\end{equation}
	for every $\varphi \in \mathcal{S}_{\sigma, \mu, C}$.
\end{proposition}

\begin{proof}
	Assume that \eqref{equivdual} is true and consider a sequence $\varphi_j \to 0$ in $\mathscr{F}_\mu$. Then, there exist $C > 0$ and $\sigma > 1$ such that $\varphi_j \in \mathcal{S}_{\sigma, \mu, C}$ for every $j$ and
	$$
	\sup_{\alpha, \beta, \gamma} \sup_{(t, x) \in \TT^m \times \R^n} C^{-|\alpha| - |\beta| - |\gamma|} (\alpha! \beta!)^{-\mu} (\gamma!)^{-\sigma} |x^\alpha \partial_x^\beta \partial_t^\gamma \varphi_j(t, x)| \to 0.
	$$
	It follows from \eqref{equivdual} that $\lim_{j \to \infty} u(\varphi_j) = 0$. \\
	\indent
	Conversely, if $u \in \mathscr{F}'_\mu$ does not satisfy \eqref{equivdual}, then there exist $C > 0$, $\sigma > 1$ and a sequence $\varphi_j \in \mathcal{S}_{\sigma, \mu, C}$ satisfying
	\begin{equation*}
		|u(\varphi_j)| > j \ \sup_{\alpha, \beta, \gamma} \ \sup_{(t, x)} C^{-|\alpha| - |\beta| - |\gamma|} (\alpha! \beta!)^{-\mu} (\gamma!)^{-\sigma} |x^\alpha \partial_x^\beta \partial_t^\gamma \varphi_j(t, x)|. \label{contrad}
	\end{equation*}
	Since $|u(\varphi_j)| > 0$, define $\psi_j \doteq \dfrac{1}{|u(\varphi_j)|}\varphi_j \in \mathcal{S}_{\sigma, \mu, C}$ for all $j \in \N$. Then
	\begin{equation*}
		\sup_{\alpha, \beta, \gamma} \ \sup_{(t, x)} C^{-|\alpha| - |\beta| - |\gamma|} (\alpha! \beta!)^{-\mu} (\gamma!)^{-\sigma} |x^\alpha \partial_x^\beta \partial_t^\gamma \psi_j(t, x)| < \frac{1}{j}, \quad j \in \N.
	\end{equation*}
	
	Therefore, $\psi_j \to 0$ in $\mathcal{S}_{\sigma, \mu, C}$, but $|u(\psi_j)| = 1$ for every $j \in \N$, which contradicts the continuity of $u$.
\end{proof}

The following result can be proved using the proposition above.

\begin{proposition}
	The restriction of any element of $\mathscr{F}'_\mu$ to $\mathcal{S}_{\sigma, \mu}$ yields an element of $\mathcal{S}'_{\sigma, \mu}$. Conversely, if $u \in \bigcap_{\sigma > 1} \mathcal{S}'_{\sigma, \mu}$, then $u \in \mathscr{F}'_\mu$.
\end{proposition}

\begin{proof}
	If $u \in \mathscr{F}'_\mu$, then its restriction on $\mathcal{S}_{\sigma, \mu}$ is linear for every $\sigma > 1$. Additionally, for each fixed $\sigma > 1$ and given $C > 0$, there exists $B_{C, \sigma} > 0$ such that \eqref{equivdual} holds for every $\varphi \in \mathcal{S}_{\sigma, \mu, C}$. This demonstrates that the restriction of $u$ on $\mathcal{S}_{\sigma, \mu}$ is an element of $\mathcal{S}'_{\sigma, \mu}$.
	\\ \indent
	Conversely, let us consider $u \in \bigcap_{\sigma>1} \mathcal{S}'_{\sigma, \mu}$. This means that $u$ is defined on all $\mathscr{F}_\mu$. Moreover, for any given $\varphi, \psi \in \mathscr{F}_\mu$, we have $\varphi, \psi \in \mathcal{S}_{\sigma, \mu}$, for some $\sigma>1$, and $u(a \varphi + b\psi) = au(\varphi) + bu(\psi)$, for any $a,b\in \C$. Therefore, $u$ is linear on $\mathscr{F}_\mu$.
	\\ \indent 
	Also, for any $\sigma>1$ and $C>0$, there exists $B_{C,\sigma}$ such that 
	$$|u(\varphi)| \leq B_{C, \sigma}\sup_{\alpha, \beta, \gamma} \sup_{(t,x) \in \TT^m \times \R^n} C^{-|\alpha|-|\beta|-|\gamma|}(\alpha! \beta!)^{-\mu} (\gamma!)^{-\sigma} |x^\alpha \partial_x^\beta \partial_t^\gamma \varphi(t,x)|,$$	
	for every $\varphi \in \mathcal{S}_{\sigma, \mu, C}.$ This implies that $u \in \mathscr{F}'_\mu$.
\end{proof}

This proposition allows us to apply to the elements of $\mathscr{F}'_{\mu}$ the characterizations in terms of the behavior of the Fourier coefficients provided in  Theorem \ref{charac_full_fourier-functions}.

\begin{lemma}\label{lemma_bound_seque}
	Let  $\{u_j\}_{j \in \N}$ be a uniformly bounded sequence in $\mathcal{G}^{\sigma}(\mathbb{T}^m)$ for some $\sigma>1$. Then, 
	$$
	u = \sum_{j \in \N} u_j(t) \varphi_j(x) \in \mathscr{F}'_{\mu}(\mathbb{T}^m \times \R^n),
	$$
\end{lemma}	
\begin{proof}
	There exists $C>0$ such that 
	$$
	\sup_{t \in \mathbb{T}^m}|u_j(t)|\leq C, \ \forall j \in \N.
	$$
	Then, given $s>1$, $h>0$ and $\psi \in \mathcal{G}^{s,h}(\mathbb{T}^m)$, we have
	$$
	|\langle u_j, \psi \rangle|  \leq 
	\int_{\mathbb{T}^m}|u_j(t) \psi(t)|dt 
	\leq (2\pi)^m C  \|\psi\|_{s, h} \leq (2\pi)^m C  \|\psi\|_{s, h} \exp\left(\varepsilon j^{\frac{1}{2n\mu}}\right),
	$$
	for all $j \in \N$ and for all $\epsilon>0$.
	Since $s>1$ is arbitrary, it follows that
	$$
	u \in \bigcap_{s>1} \mathcal{S}_{s, \mu}'(\mathbb{T}^m \times \R^n) = \mathscr{F}'_{\mu}(\mathbb{T}^m \times \R^n).
	$$
\end{proof}

\section{Reduction to the normal form}\label{secRedNormalForm}

To prove Theorem \ref{main_Theorem}, we first reduce the system to a simpler form. Specifically, we will show that the global hypoellipticity of the original system is equivalent to that of a reduced system, referred to as the normal form.
In the special case where the time dependence involves only real-valued functions, the normal form simplifies to a time-independent system, a scenario treated in \cite{AviCapKiri24}.

\subsection{Automorphisms and hypoellipticity}

Consider the system $\LL = (L_1, L_2, \ldots, L_m)$ defined on $\mathbb{T}^m \times \mathbb{R}^n$, given by
\begin{equation}\label{L_a}
	L_r = D_{t_r} + (a_r+ib_r)(t_r)P(x,D_x), \quad r=1, \ldots, m,
\end{equation}
where $a_{r}, b_r \in \mathcal{G}^{\sigma}(\mathbb{T}^1; \mathbb{R})$ are real-valued functions, the derivatives $D_r = -i \partial/\partial t_r$ are defined on distinct copies of the one-dimensional torus, and the operator $P = P(x, D_x)$ satisfies the conditions \eqref{P-intro1} and \eqref{P-elliptic-intro}.

The normal form of $\LL$ is the system $\LL_0 = (L_{1,0}, L_{2,0}, \ldots, L_{m,0})$ defined on $\mathbb{T}^m \times \mathbb{R}^n$, given by
\begin{equation}\label{L_0}
	L_{r,0} = D_{t_r} + (a_{r,0} + ib_r(t_r))P(x,D_x), \quad r = 1, \ldots, m,
\end{equation}
where
$a_{r,0}$ is defined by \eqref{realpartaverage}.

This method of deducing the properties of $\LL$ from those of its normal form $\LL_{0}$ is widely used for systems and operators on the torus and compact Lie groups, as noted in \cite{AriKirMed19, AviGraKir18, Ber99, KirMorRuz20, KirMorRuz21jfa, Pet05}. The idea is to find a linear automorphism $\Psi$ on $\mathscr{F}_\mu$ and on $\mathscr{F}'_\mu$ such that
\begin{equation}\label{conjugationformula}
	\Psi^{-1} \circ L_r \circ \Psi = L_{r,0}, \quad 1 \leq r \leq m.
\end{equation}

Assuming the existence of such an automorphism, if $\LL_{0}$ is $\mathscr{F}_{\mu}$-globally hypoelliptic and $u \in \mathscr{F}_{\mu}'$  is a solution of $L_{r} u = f_r \in \mathscr{F}_{\mu}$ for each $r = 1, \ldots, m$, then writing 
$$
f_r = \Psi g_r, \quad \text{with } g_r \in \mathscr{F}_{\mu} \text{ for all } r = 1, \ldots, m.
$$
by \eqref{conjugationformula}, we have
$$
L_{r} u = \Psi g_r \Longrightarrow L_{r,0}(\Psi^{-1} u) = g_r, \quad r = 1, \ldots, m.
$$

Since $\Psi^{-1} u \in \mathscr{F}_{\mu}'$, it follows from the hypothesis on $\LL_{0}$ that $\Psi^{-1} u \in \mathscr{F}_{\mu}$, which implies that $u \in \mathscr{F}_{\mu}$. Therefore, $\LL$ is $\mathscr{F}_{\mu}$-globally hypoelliptic. Similarly, the $\mathscr{F}_{\mu}$-global hypoellipticity of $\LL$ implies the same property for $\LL_{0}$, using a similar argument. This argument leads to the following result:

\begin{theorem}\label{The-Normal}
	Let $\sigma > 1$ and $\mu \geq 1/2$. Then the system $\mathbb{L}$ is $\mathscr{F}_{\mu}$-globally hypoelliptic if and only if its normal form $\mathbb{L}_{0}$ is $\mathscr{F}_{\mu}$-globally hypoelliptic.
\end{theorem}

To complete the proof of Theorem \ref{The-Normal}, it remains to establish the existence of an automorphism that satisfies \eqref{conjugationformula}, which is provided in the following proposition.

\begin{proposition}\label{Theorem-Psi}
	Let $u(t, x) = \sum_{j \in \N} u_j(t) \varphi_j(x) \in \mathcal{S}_{\sigma, \mu}'$, and define
	\begin{equation}\label{conjugation}
		(\Psi u)(t, x) \doteq \sum_{j \in \N} \exp\left(-i A(t) \lambda_j\right) u_j(t) \varphi_j(x),
	\end{equation}
	where
	\begin{equation}\label{A(t)definition}
		A(t) = \sum_{k = 1}^{m} \left(\int_{0}^{t_k} a_k(\eta) \, d\eta - a_{k, 0} t_k \right).
	\end{equation}
	Then, $\Psi: \mathscr{F}_{\mu} \to \mathscr{F}_{\mu} $ and $\Psi: \mathscr{F}_{\mu}' \to \mathscr{F}_{\mu}'$
	are automorphisms, with the inverse given by
	\begin{equation}\label{conjugation_inverse}
		\Psi^{-1} u(t, x) = \sum_{j \in \N} \exp\left(i A(t) \lambda_j\right) u_j(t) \varphi_j(x).
	\end{equation}
	Moreover, \eqref{conjugationformula} holds.
\end{proposition}

\begin{proof}
	We first show that $\Psi: \mathscr{F}_{\mu} \to \mathscr{F}_{\mu}$ is well-defined. Given $f = \sum_{j \in \mathbb{N}} f_j(t) \varphi_j(x) \in \mathcal{S}_{\delta, \mu}$, denote $\psi_j(t) = \exp\left(-i A(t) \lambda_j\right) f_j(t)$ for $j \in \N$.
	
	Then, applying the Leibniz formula for any multi-index $\alpha \in \N^m_0$, we have
	$$
	|\partial_t^{\alpha} \psi_j(t)| \leq \sum_{\beta \leq \alpha}
	\binom{\alpha}{\beta} 
	|\partial_t^{\beta} \exp \left(-i A(t) \lambda_j\right)|
	|\partial_t^{\alpha - \beta} f_j(t)|.
	$$
	
	Since $f \in \mathcal{S}_{\delta, \mu}$, there exist constants $\epsilon, C > 0$ such that
	\begin{equation*}
		\sup_{t \in \mathbb{T}^m} | \partial_t^{\alpha-\beta} f_j(t)| \leq
		C^{|\alpha-\beta|+1} [(\alpha-\beta)!]^{\delta} \exp \left[-\epsilon j^{\frac{1}{2n\mu}} \right],
	\end{equation*}    
	where $\beta = (\beta_1, \ldots, \beta_m) \in \N_0^m$. Note that
	$$
	|\partial_t^{\beta} \exp \left(-i A(t) \lambda_j\right)| = 
	\prod_{k=1}^{m} \left| \partial_{t_k}^{\beta_k} \exp \left[-i \lambda_j \left(\int_{0}^{t_k} a_k(\eta) \, d\eta - a_{k, 0} t_k \right)\right]\right|.
	$$
	Denoting
	$$
	\mathcal{H}_j(t_k) = \exp\left[-i \lambda_j \left( \int_{0}^{t_k} a_k(\eta) \, d\eta - a_{k, 0} t_k \right)\right],
	$$
	we obtain from Fa\`a di Bruno's formula
	$$
	\partial_{t_k}^{\beta_k} \mathcal{H}_j(t_k) = \sum_{\Delta(\gamma), \, \beta_k}
	\frac{1}{\gamma!}(-i \lambda_j)^\gamma
	\frac{\beta_k!}{s_1! \cdots s_\gamma!}
	\left( \prod_{\nu=1}^\gamma
	\partial_{t_k}^{s_\nu - 1}(a_k(t_k) - a_{k, 0}) \right)
	\mathcal{H}_j(t_k),
	$$
	where the summation over $\Delta(\gamma)$ and $\beta_k$ ranges over all integers $1 \leq \gamma \leq \beta_k$, and the sum is taken over all $s_\nu \in \N$ such that $s_1 + \ldots + s_\gamma = \beta_k$.
	
	Since
	$$
	\left| \prod_{\nu=1}^\gamma
	\partial_{t_k}^{s_\nu-1}(a_k(t_k)-a_{k,0})  \right| \leq C^{\beta_k-\gamma+1}[(\beta_k -\gamma)!]^\sigma,
	$$
	and $|\lambda_j| \leq C_1 j^{M /2n}$, we get
	\begin{align*}
		|\partial_{t_k}^{\beta_k} \mathcal{H}_j(t_k)| & \leq 
		\sum_{\Delta(\gamma), \, \beta_k}|\lambda_j|^\gamma
		\dfrac{1}{\gamma!}
		\frac{\beta_k!}{s_1! \cdots s_\gamma! }
		C^{\beta_k-\gamma+1}[(\beta_k -\gamma)!]^\sigma \\
		& \leq 
		\sum_{\Delta(\gamma), \, \beta_k}C_1^{\gamma}j^{\frac{\gamma M}{2n}}
		\dfrac{1}{\gamma!} 
		\frac{\beta_k!}{s_1! \cdots s_\gamma! }
		C^{\beta_k-\gamma+1}[(\beta_k -\gamma)!]^\sigma 
	\end{align*}
	
	In view of the inequality
	\begin{equation*}
		j^{\frac{\gamma M}{2n}} 
		\leq C_{\epsilon}^{\gamma}(\gamma !)^{M\mu} \exp \left(\epsilon/(2m) j^{\frac{1}{2n\mu}} \right),
	\end{equation*}
	we conclude that
	$$
	|\partial_{t_k}^{\beta_k} \mathcal{H}_j(t_k)| \leq C_2^{\beta_k + 1} \exp \left(\epsilon/(2m)j^{\frac{1}{2n\mu}} \right) (\beta_k !)^{\max\{\sigma, M \mu -1\}}.
	$$ 
	
	Therefore, 
	\begin{equation}\label{est_e_A}
		|\partial_t^{\beta}\exp \left(-i A(t) \lambda_j\right)| \leq 
		C_2^{|\beta|+ 1}
		\exp \left(\epsilon/2 \, j^{\frac{1}{2n\mu}} \right) (\beta!)^{\max\{\sigma, M \mu -1\}}
	\end{equation}
	implying
	$$
	|\partial^{\alpha}_t \psi_j(t)|\leq C^{|\alpha| +1} (\alpha !)^{\rho} 
	\exp \left(-\epsilon/2 \, j^{\frac{1}{2n\mu}} \right),
	$$	 
	where 
	$$
	\rho = \max\{\delta, \max\{\sigma, M \mu -1\}\}
	$$
	and consequently, $\Psi f \in \mathscr{F}_{\mu}$.
	
	\medskip
	Now, let us show that if $u \in \mathscr{F}_{\mu}'$, then $\Psi u \in \mathscr{F}_{\mu}'$. To do this, denoted $\psi_j(t)= \exp(-iA(t)\lambda_j)u_j(t)$, it is sufficient to prove that, given $\sigma > 0$, we can find a constant $C_{\epsilon, h} > 0$ for every $\epsilon, h > 0$ satisfying
	\begin{equation*}
		|\langle  \psi_j(t) \, , \, \theta(t) \rangle | \leq C_{\epsilon,h} \|\theta\|_{\sigma, h} \exp \left(\epsilon j^{\frac{1}{2n\mu}}\right),
	\end{equation*}
	for all $j \in \mathbb{N}$, and for all $\theta \in \mathcal{G}^{\sigma,h}(\TT^m)$, cf. \cite[Theorem 2.9]{AviCap22}.
	
	Let  $u \in \mathscr{F}_{\mu}'$ be fixed. In view of  \eqref{secondnorm}, we have the following: for every $\eta>0$ and $A>0$, there is $B_{A,\eta}>0$ such that
	\begin{equation*}
		|\langle  u\, , \, \Theta(t,x) \rangle | 
		\leq  B_{A,\eta} \sup\limits_{\substack{\alpha \in \N^m \\ k \in \N}}    
		A^{-k - |\alpha|} (k!)^{-M \mu} (\alpha!)^{-\eta} \|P^k \partial_t^{\alpha}\Theta\|_{L^2(\TT^m\times \R^n)},
	\end{equation*}
	for all $\Theta \in \mathcal{S}_{\sigma,\mu,C}$.

	Let $\epsilon, h>0$ be fixed. Given $\theta \in \mathcal{G}^{\sigma,h}(\TT^m)$, it follows from \eqref{est_e_A} that
	\begin{align*}
		& |\partial_t^{\alpha}\left[\exp\left(-i A(t) \lambda_j\right)\theta(t)\right]|\\[2mm]
		&\ \leq   \exp\left(\epsilon/2 j^{\frac{1}{2n\mu}} \right)
		\sum_{\beta \leq \alpha} \binom{\alpha}{\beta}  C^{|\beta|+ 1}
		(\beta!)^{\max\{\sigma, M\mu -1\}} |\partial_t^{\alpha -\beta}\theta(t)| \\
		&\ \leq  \exp\left(\epsilon/2 j^{\frac{1}{2n\mu}} \right) \sum_{\beta \leq \alpha} \binom{\alpha}{\beta}  C^{|\beta|+ 1}
		(\beta!)^{\max\{\sigma, M\mu -1\}} \|\theta\|_{\sigma,h}((\alpha - \beta)!)^{\sigma}h^{|\alpha-\beta|} \\
		&\ \leq  \exp\left(\epsilon/2 j^{\frac{1}{2n\mu}} \right) C_{\epsilon,h}^{|\alpha|} (\alpha!)^{\rho} \|\theta\|_{\sigma,h},
	\end{align*} 
	where $\rho = \max\{\sigma,  M\mu -1 \}$.
	
	\medskip
	Hence, setting
	$$
	\Omega_{j}(t)=\langle  \psi_j(t) \, , \, \theta(t) \rangle = \langle  u \, , \, \exp\left(-i\lambda_j A(t) \right)\theta(t)\varphi_j(x) \rangle,
	$$
	we obtain
	\begin{align*}
		|\Omega_{j}(t)| 
		& \leq 
		B_{A,\eta}  \sup\limits_{\stackrel{\alpha \in \N^m_0}{k \in \N}}	
		A^{-k - |\alpha|} (k!)^{-M \mu} (\alpha!)^{-\eta}\cdot \|P^k \partial_t^{\alpha}\left[\exp\left(-i A(t) \lambda_j\right)\theta(t)\varphi_j(x) \right]\|_{L^2(\TT^m\times \R^n)}\\
		& = 
		B_{A,\eta} \sup\limits_{\stackrel{\alpha \in \N^m_0}{k \in \N}}	
		A^{-k - |\alpha|} (k!)^{-M \mu} (\alpha!)^{-\eta}  |\partial_t^{\alpha}\left[\exp\left(-i A(t) \lambda_j\right)\theta(t)\right]| \cdot \|P^k\varphi_j(x) \|_{L^2(\TT^m\times \R^n)} \\
		& =
		B_{A,\eta} \sup\limits_{\stackrel{\alpha \in \N^m_0}{k \in \N}}	
		A^{-k - |\alpha|} (k!)^{-M \mu} (\alpha!)^{-\eta}  |\partial_t^{\alpha}\left[\exp\left(-i A(t) \lambda_j\right)\theta(t)\right]| \cdot |\lambda_j|^k \\
		& \leq
		B_{A,\eta} \sup\limits_{\stackrel{\alpha \in \N^m_0}{k \in \N}}	
		A^{-k - |\alpha|} (k!)^{-M \mu} (\alpha!)^{-\eta}   C^{k}j^{kM/2n } |\partial_t^{\alpha}\left[\exp\left(-i A(t) \lambda_j\right)\theta(t)\right]|  \\
		& \leq
		B_{A,\eta} \sup\limits_{\stackrel{\alpha \in \N^m_0}{k \in \N}}	
		A^{-k - |\alpha|} (k!)^{-M \mu} (\alpha!)^{-\eta}   C^{k} C_{\epsilon}^{k}(k!)^{M\mu} \exp\left(\epsilon j^{\frac{1}{2n\mu}} \right) C_{\epsilon,h}^{|\alpha|} (\alpha!)^{\rho} \|\theta\|_{\sigma,h}  \\
		& = 
		B_{A,\eta} \|\theta\|_{\sigma,h} \exp\left(\epsilon j^{\frac{1}{2n\mu}} \right) \sup\limits_{\alpha \in \N^m_0}\sup\limits_{k \in \N} A^{-k - |\alpha|} (\alpha!)^{\rho-\eta}  C_{\epsilon, h}^{|\alpha|} C_{\epsilon}^{k}.
	\end{align*}
	
	Therefore, for $\eta \geq \rho$ and 
	$A^{-1} = \max\{C_{\epsilon}, C_{\epsilon,h}\}$,  
	$$
	|\langle  \psi_j(t) \, , \, \theta(t) \rangle|\leq \widetilde{C}_{\epsilon,h}  \|\theta\|_{\sigma,h} \exp\left(\epsilon j^{\frac{1}{2n\mu}} \right).
	$$
	
	To complete the proof, it remains to verify that $\Psi^{-1}$ is indeed the inverse of $\Psi$ and that both operators are linear. Finally, from equations \eqref{conjugation} and \eqref{conjugation_inverse}, it follows clearly that $L_r \circ \Psi = \Psi \circ L_{r,0}$ for every $1 \leq r \leq m$.
\end{proof}

\subsection{Time independent coefficients}
Consider the system with time-independent coefficients
$$
\mathscr{L} = \{\mathscr{L}_r = D_{t_r} + \omega_{r} P(x, D_x), \ r= 1,\ldots, m\}
$$	
where $\omega_{r} \in \C$. Also, we set
$$
\sigma_{\mathscr{L}}(\tau,j) =  (\sigma_{\mathscr{L}_{1}}(\tau_1,j), \ldots, \sigma_{\mathscr{L}_{m}}(\tau_m,j)), \ \tau = (\tau_1, \ldots, \tau_m) \in \Z^m, \ j \in \N, 
$$
where $\sigma_{\mathscr{L}_{r}}(\tau_r,j) = \tau_r + \omega_r \lambda_j.$

We also define
$$
\|\sigma_{\mathscr{L}}(\tau,j)\| \doteq \max_{1 \leq r \leq m} |\sigma_{\mathscr{L}_r}(\tau,j)|, \quad (\tau,j) \in \mathbb{Z}^m \times \mathbb{N},
$$
and the set
\begin{equation*}
	\mathcal{N}_{\mathscr{L}} = \{(\tau,j) \in \mathbb{Z}^N \times \mathbb{N}; \, \sigma_{\mathscr{L}}(\tau,j) = 0 \}.
\end{equation*}

\begin{theorem} \label{GH-thm}
	The system $\mathscr{L}$ is $\mathscr{F}_{\mu}$-globally hypoelliptic if and only if the set $\mathcal{N}_{\mathscr{L}}$ is finite and there is  $\sigma_0>1$ such that, for all $\sigma\geq \sigma_0$ we have the following property: for all $\varepsilon > 0$, there exists $C_\varepsilon > 0$ such that  
	\begin{equation}\label{Dio-cond-time-ind}
		\|\sigma_{\mathscr{L}}(\tau,j)\| \geq C_\varepsilon \exp\left[-\varepsilon(\|\tau\|^{1/\sigma} + j^{1/(2n\mu)})\right],
	\end{equation}
	for any $(\tau,j) \in \mathbb{Z}^m \times \mathbb{N}$ such that $\sigma_{\mathscr{L}}(\tau,j) \neq 0$.
\end{theorem}

\begin{proof} We proceed with an argument similar to that in the proof of \cite[Theorem 1.2]{AviCapKiri24}. We start with the sufficiency. Assume that $\mathcal{N}_{\mathscr{L}}$ is finite and that \eqref{Dio-cond-time-ind} holds. Let $f_r \in \mathscr{F}_{\mu}$, $r=1, \ldots, m$ such that
	$$
	\mathscr{L}_{r} u = f_r, \textrm{ for some } u \in \mathscr{F}_{\mu}'.
	$$ 
	We may assume that there exists $\widetilde{\sigma} >1$ such that $f_r \in \mathcal{S}_{\widetilde{\sigma}, \mu}$ for all $r=1,\ldots, m$. For each  
	$(\tau,j) \in \mathbb{Z}^m\times \mathbb{N}$ we set $r^* = r^*(\tau,j)\in\{1,\ldots,m\}$ as the index such that 
	$$
	|\sigma_{\mathscr{L}_{r^*}}(\tau,j)| = \|\sigma_{\mathscr{L}}(\tau,j)\| = \max_{1\leq r\leq m} |\sigma_{\mathscr{L}_r}(\tau,j)|.
	$$

	Due to the finiteness of the set $\mathcal{N}_{\mathscr{L}}$, it is sufficient to analyze the behavior of $\widehat{u_j}(\tau)$ outside $\mathcal{N}_{\mathscr{L}}$. If $(\tau, j) \notin \mathcal{N}_{\mathscr{L}}$, then $\widehat{u_j}(\tau) = \sigma_{\mathscr{L}_{r^*}}(\tau,j)^{-1} \widehat{f_{r^*,j}}(\tau)$, leading to
	\begin{align*}
		|\widehat{u_j}(\tau)| &= |\sigma_{\mathscr{L}_{r^*}}(\tau,j)|^{-1}|\widehat{f_{r^*,j}}(\tau)| \\ 
		&= \|\sigma_{\mathscr{L}}(\tau,j)\|^{-1} |\widehat{f_{r^*,j}}(\tau)| \\
		&\leq C_\varepsilon  \exp\left[\varepsilon (\|\tau\|^{1/\sigma}+j^{1/(2n\mu)}) \right]|\widehat{f_{r^*,j}}(\tau)|,
	\end{align*}
	where $\sigma\geq \widetilde{\sigma}$.
	
	Since $f_{r^*} \in \mathcal{S}_{\widetilde{\sigma},\mu}$, then there is   $\varepsilon_{r^*}>0$ and $C_{r^*}>0$ such that 
	$$ |\widehat{f_{{r^*},j}}(\tau)| \leq C_{r^*} \exp\left[-\varepsilon_{r^*}(\|\tau\|^{1/\widetilde{\sigma}}+j^{1/(2n\mu)})\right],
	$$    
	for every $(\tau,j)\in\Z^m\times \N.$
	
	Therefore, for $\varepsilon=\varepsilon_{r^*}/2$ we have
	$$
	|\widehat{u_j}(\tau)|  \leq C_\varepsilon  C_{r^*} \exp\left[-\frac{\varepsilon_{r^*}}{2}(\|\tau\|^{1/\widetilde{\sigma}}+j^{1/(2n\mu)})\right].
	$$
	
	Thus, $u \in \mathcal{S}_{\widetilde{\sigma},\mu}$, implying that $\mathscr{L}$ is $\mathscr{F}_{\mu}$-globally hypoelliptic.

	For the converse we preliminarly observe that if $\mathcal{N}_{\mathscr{L}}$ is infinite one can easily construct a solution $u \in \mathscr{F}'_\mu \setminus \mathscr{F}_\mu$ of the system $\mathscr{L} u=0$ by defining 
	\begin{equation}\label{proof-N-infinite}
		\widehat{u}_j(\tau)= \begin{cases} 1, \quad \textrm{if} \quad \sigma_{\mathscr{L}}(\tau,j)=0 \\ 0, \quad \textrm{otherwise} \end{cases},
	\end{equation}
	and conclude that $\mathscr{L}$ is not $\mathscr{F}_\mu$-globally hypoelliptic. To conclude the proof, let $\mathscr{L}$ be $\mathscr{F}_\mu$-globally hypoelliptic and assume the existence of $\sigma >1$ and 
	$\varepsilon >0$ such that for each $k\in\N$ there exists $(\tau_k,j_k)\in\Z^m\times\N$ such that
	\begin{equation*}
		0<\|\sigma_{\mathscr{L}}(\tau_k,j_k)\| <  \exp\left[-\varepsilon (\|\tau_k\|^{1/\sigma} + j_k^{1/(2n\mu)}) \right].
	\end{equation*}
	
	Now, for each $r=1,\ldots,\ell$, we define 
	$$
	\widehat{f_{r,j}}(\tau)= \begin{cases}
		\sigma_{\mathscr{L}_r}(\tau_k,j_k) , & \text{if } (\tau,j) = (\tau_k,j_k) \text{ for some } k\in \N, \\
		0, & \text{otherwise. }
	\end{cases}
	$$

	Thus $f_r =  \sum\limits_{j \in \N,\tau \in\Z^m} \widehat{f_{r,j}}(\tau)e^{i\tau \cdot t}\varphi_{j}(x) \in \mathcal{S}_{\sigma,\mu}$ for every $r=1,\ldots,m$.  On the other hand, by defining  
	$$
	\widehat{u_j}(\tau) = 
	\begin{cases}
		1 , & \text{if } (\tau,j) = (\tau_k,j_k) \text{ for some } k\in \N, \\
		0, & \text{otherwise. }
	\end{cases}
	$$
	we get
	$$
	u(t,x)= \sum_{j \in \N}\sum_{\tau \in\Z^m} \widehat{u_j}(\tau)e^{i\tau \cdot t}\varphi_{j}(x) \in \mathcal{S}_{\sigma,\mu}' \setminus \mathcal{S}_{\sigma,\mu}
	$$
	for all $\sigma>1$. Therefore, by Theorem \ref{charac_full_fourier-functions}, $u\in \mathscr{F}_ {\mu}' \setminus \mathscr{F}_{\mu}$ and  $\mathscr{L}_{r} u=f_r$, $r=1, \ldots, m$, implying that $\mathscr{L}$ is not $\mathscr{F}_{\mu}$-globally hypoelliptic.
\end{proof}

We now specialize our attention to the case where the coefficients are real-valued and time-dependent, aiming to characterize the global hypoellipticity for such systems. 

Consider the system $\LL = (L_1, L_2, \ldots, L_m)$ defined on $\mathbb{T}^m \times \mathbb{R}^n$, given by
\begin{equation}\label{real-case-variable}
	L_r = D_{t_r} + a_r(t_r)P(x,D_x), \quad r=1, \ldots, m,
\end{equation}
where each $a_{r} \in \mathcal{G}^{\sigma}(\mathbb{T}^1; \mathbb{R})$ is a real-valued functions. Its normal form, as defined previously (see \eqref{L_0} with $b_r \equiv 0$), is the constant-coefficient system \( \mathbb{L}_0 = (L_{1,0}, L_{2,0}, \ldots, L_{m,0})\) with
\begin{equation}\label{real-case-normal-form}
	L_{r,0} = D_{t_r} + a_{r,0} P(x, D_x), \ r= 1,\ldots, m,
\end{equation}
where
$ a_{r,0}$ is given by \eqref{realpartaverage}.

Combining Theorems~\ref{The-Normal} and \ref{GH-thm} yields the following result for this class of systems:

\begin{proposition}\label{Prop_real_valued_coeff}
	The system \eqref{real-case-variable} is $\mathscr{F}_{\mu}$-globally hypoelliptic if and only if its normal form, the system \eqref{real-case-normal-form} is $\mathscr{F}_{\mu}$-globally hypoelliptic,This, in turn, is equivalent to the conditions that the set $\mathcal{N}_{\mathbb{L}_0}$ of zeros of $\LL_{0}$ is finite, and that there is  $\sigma_0>1$ such that for all $\sigma\geq \sigma_0$ the following condition holds: for every $\varepsilon > 0$, there exists $C_\varepsilon > 0$ such that  
	\begin{equation*}
		\|\sigma_{\LL_0}(\tau,j)\| \geq C_\varepsilon \exp\left[-\varepsilon(\|\tau\|^{1/\sigma} + j^{1/(2n\mu)})\right],
	\end{equation*}
	for all $(\tau,j) \in \mathbb{Z}^m \times \mathbb{N}$ such that $\sigma_{\LL_0}(\tau,j) \neq 0$.
\end{proposition}

\section{Sufficient conditions for the hypoellipticity}\label{secSuffCOnditions}
In this section we prove the sufficiency of conditions (I) and (II) in Theorem \ref{main_Theorem} for the $\mathscr{F}_\mu$-global hypoellipticity of the system \eqref{general-const-system}.
\\ \indent
Let us start with some preliminary results. Let $c_{r,0} = (2\pi)^{-1}\int_0^{2\pi} c_r(t_r)\, dt_r, r=1,\ldots, m,$ and consider the sets
\begin{equation}\label{Z_L-set}
	\mathcal{Z}_r = \{j \in \N \mid c_{r,0} \lambda_j \in \Z\},  \text{ and } \
	\mathcal{Z}_{\LL} = \bigcap_{r=1}^{m} \mathcal{Z}_r.
\end{equation}

\begin{proposition}\label{Z-infinity}
	If $\mathcal{Z}_{\LL}$ is an infinite set, then $\LL$ is not $\mathscr{F}_\mu$-globally hypoelliptic.
\end{proposition}

\begin{proof}
	Let $\{j_{\ell}\}_{\ell \in \N}$ be an increasing sequence in $\mathcal{Z}_{\LL}$, and assume $\lambda_{j_\ell} > 0$. For each $r = 1, \ldots, m$, let $t^*_r \in [0, 2\pi]$ be the value such that
	$$
	\int_{0}^{t^*_r} b_{r}(\zeta) \, d\zeta =
	\max_{0 \leq t \leq 2\pi} \int_{0}^{t} b_{r}(\zeta) \, d\zeta.
	$$
	For each $\ell \in \N$, consider the periodic function
	$$
	u_{\ell}(t) = \prod_{r = 1}^{m} \exp\left[-\lambda_{j_{\ell}} \left(i \int_{0}^{t_r} c_{r}(\zeta) \, d\zeta
	+ \int_{0}^{t^*_r} b_{r}(\zeta) \, d\zeta \right)\right].
	$$
	
	Observe that $|u_\ell(t^*_1, \ldots, t^*_m)| = 1$ and $|u_{\ell}(t)| \leq 1$ for all $\ell \in \N$ and $t \in \TT^m$. With this, it is not difficult to prove that
	\begin{equation}\label{solution-Z-infinity}
		u(t,x) = \sum_{\ell \in \N} u_{\ell}(t) \varphi_\ell(x)  \in \mathscr{F}_\mu' \setminus \mathscr{F}_\mu.
	\end{equation}
	
	Moreover, for all $\ell \in \N$ and $t \in \TT^m$, we have
	$$
	D_{t_r} u_{\ell}(t) + \lambda_{j_\ell} c_{r}(t_r) u_{\ell}(t) = 0, \quad r = 1, \ldots, m,
	$$
	which implies that $\LL u = 0$ and $\LL$ is not $\mathscr{F}_\mu$-globally hypoelliptic.
\end{proof}

\begin{remark}\label{lambda-I}
	We note that with a slight modification of the argument above, we can address the general case where  $\lambda_{j}$ is not necessarily positive. Indeed, consider the sets
	\begin{equation}\label{setsW}
		\mathcal{W}_{+} = \{j \in \N \mid \lambda_{j} > 0\} \quad \text{and} \quad
		\mathcal{W}_{-} = \{j \in \N \mid \lambda_{j} < 0\}.
	\end{equation}
	
	The proof of Proposition \ref{Z-infinity} covers the case where $\mathcal{W}_{+}$ is an infinite set. If $\mathcal{W}_{+}$ is finite, then $\mathcal{W}_{-}$ is infinite and it is sufficient to define
	$$
	u_{\ell}(t) = \prod_{r = 1}^{m} \exp\left[-\lambda_{j_{\ell}} \left(i \int_{0}^{t_r} c_{r}(\zeta) \, d\zeta
	+ \int_{0}^{\tilde{t}_r^*} b_{r}(\zeta) \, d\zeta \right)\right],
	$$
	where $\tilde{t}_r^*$ is such that
	$$
	\int_{0}^{\tilde{t}_r^*} b_{r}(\zeta) \, d\zeta =
	\min_{0 \leq t \leq 2\pi} \int_{0}^{t} b_{r}(\zeta) \, d\zeta.
	$$
	and repeat the same steps as in the previous proof.
\end{remark}

Now, taking into account Remark \ref{GHcomparison}, we reformulate the result obtained in \cite{AviCap22} for the case $m=1$ in terms of the notion of global hypoellipticity adopted in this paper.

\begin{proposition}\label{GHm=1}
	The operator 
	$$\mathcal{L} = D_t + (a+ib)(t)P(x, D_x)$$ 
	is $\mathscr{F}_{\mu}$-globally hypoelliptic if and only if one of the following conditions holds:
	\begin{enumerate}
		\item[$(a)$] $b$ is not identically zero and does not change sign;
		\item[$(b)$] $b \equiv 0$ and $a_0 = (2\pi)^{-1} \int_{0}^{2\pi} a(s) \, ds$ satisfies $\mathscr{D}_{\sigma,\mu}$ for all $\sigma\geq M\mu$.
	\end{enumerate}
	
\end{proposition}

To prove the proposition we need the following lemma.


\begin{lemma}\label{lt2}
	Let $\mu \geq \frac{1}{2}$, $\sigma \geq M\mu$, and $\omega = (\omega_1, \ldots, \omega_m) \in \mathbb{R}^m \setminus \{0\}$. The following statements are equivalent:
	\begin{itemize}
		\item[(i)] For every $\epsilon > 0$, there exists a constant $C_\epsilon > 0$ such that
		\[
		\max_{1 \leq r \leq m} |\tau_r - \omega_r \lambda_j| \geq C_\epsilon \exp\left\{-\epsilon\left( \|\tau\|^{\frac{1}{\sigma}} + j^{\frac{1}{2n\mu}} \right) \right\}, \quad \forall \tau \in \mathbb{Z}^m,\ \forall j \in \mathbb{N}.
		\]
		
		\item[(ii)] For every $\delta > 0$, there exists a constant $C_\delta > 0$ such that
		\[
		\max_{1 \leq r \leq m} |1 - \exp (-2\pi i \omega_r \lambda_j)| \geq C_\delta \exp\left\{ -\delta j^{\frac{1}{2n\mu}} \right\}, \quad \forall j \in \mathbb{N}.
		\]
	\end{itemize}
\end{lemma}

\begin{proof}
	Suppose (i) fails. Then there exist $\epsilon > 0$, a sequence $(\tau_k, j_k) \in \mathbb{Z}^m \times \mathbb{N}$, such that
	\begin{align*}
		0 < \max_{1 \leq r \leq m} |\tau_{k,r} - \omega_r \lambda_{j_k}| &< \exp\left\{-\epsilon\left( \|\tau_k\|^{1/\sigma} + j_k^{1/(2n\mu)} \right) \right\} \\
		&\leq \exp\left\{-\epsilon j_k^{1/(2n\mu)}\right\}.
	\end{align*}
	
	Now, recall that the function
	\[
	\zeta \mapsto \frac{1 - e^{-i\zeta}}{\zeta}
	\]
	has a continuous extension near $\zeta = 0$. Hence, there exists a constant $C > 0$ such that
	\[
	|1 - e^{-i\zeta}| \leq C |\zeta|, \quad \text{for } |\zeta| < C^{-1}.
	\]
	
	Thus, for each $r = 1, \ldots, m$ and sufficiently large $k$, we have:
	\begin{align*}
		|1 - \exp{(-2\pi i \omega_r \lambda_{j_k}})| 
		&= |1 - \exp{(-2\pi i (\tau_{k,r} - \omega_r \lambda_{j_k})})| \\
		&\leq 2\pi C |\tau_{k,r} - \omega_r \lambda_{j_k}| \\
		&\leq 2\pi C \exp\left(-\epsilon j_k^{1/(2n\mu)}\right),
	\end{align*}
	contradicting (ii).

	Conversely, suppose (ii) fails. Then there exists $\epsilon > 0$ such that for any $C > 0$, there is $j_C \in \mathbb{N}$ satisfying
	$$
	0 < \max_{1 \leq r \leq m} |1 - \exp(-2\pi i \omega_r\lambda_{j_C})| < C \exp\left(-\epsilon j_C^{1/(2n\mu)}\right).
	$$
	
	Choosing $C_k = \exp\left(-\epsilon k^{1/\sigma}\right)$ for each $k \in \mathbb{N}$, we obtain a sequence of elements $j_k \in \mathbb{N}$ such that
	$$
	0 < \max_{1 \leq r \leq m} |1 - \exp(-2\pi i \omega_r\lambda_{j_k})| < \exp\left\{-\epsilon\left(k^{1/\sigma} + j_k^{1/(2n\mu)}\right)\right\}.
	$$
	In particular, this implies that $\exp(-2\pi i \omega_r\lambda_{j_k}) \to 1$ as $k \to \infty$, for all $r = 1, \ldots, m$.
	
	We now claim that there exist sequences of elements $\gamma_{r,k} \in \mathbb{Z}$, for each $r = 1, \ldots, m$, such that $|\gamma_{r,k}| \to \infty$ and
	$$
	|\gamma_{r,k} - \omega_r\lambda_{j_k}| \to 0 \quad \text{as } k \to \infty.
	$$
	Indeed, define $\gamma_{r,k} \in \mathbb{Z}$ by
	$$
	|\gamma_{r,k} - \omega_r\lambda_{j_k}| = \mathrm{dist}(\omega_r\lambda_{j_k}, \mathbb{Z}).
	$$
	
	Suppose that $ |\gamma_{r,k} - \omega_r\lambda_{j_k}| \not\to 0 $. Then there exists a constant $ c > 0 $ such that, for infinitely many $ k $, we have
	$$
	c < |\gamma_{r,k} - \omega_r\lambda_{j_k}| < \frac{1}{2}.
	$$
	
	Thus  $2\pi c < |2\pi \gamma_{r,k} - 2\pi \omega_r\lambda_{j_k}| < \pi,$ for infinitely many $ k $. This implies that the argument of the exponential function $ \exp(-2\pi i \omega_r\lambda_{j_k}) $ stays uniformly far from an integer multiple of $ 2\pi $, so the sequence cannot converge to 1. This contradicts our earlier conclusion that
	$$
	\exp(-2\pi i \omega_r\lambda_{j_k}) \to 1 \quad \text{as } k \to \infty.
	$$
	Hence, we must have
	$$
	|\gamma_{r,k} - \omega_r\lambda_{j_k}| \to 0 \quad \text{as } k \to \infty.
	$$

	Now, consider the sequence of vectors 
	$\gamma_k = (\gamma_{1,k}, \ldots, \gamma_{m,k} ) \in \Z^m.$ 
	Since the function
	$$
	\zeta \mapsto \frac{\zeta}{1 - \exp(-i\zeta)}
	$$
	has a continuous extension near $\zeta = 0$, there exists a constant $\delta > 0$ such that
	$$
	|\zeta| \leq \delta |1 - \exp(-i\zeta)|, \quad \text{for } |\zeta| < \delta^{-1}.
	$$
	
	For sufficiently large $k$, we have $2\pi |\gamma_{r,k} - \omega_{r}\lambda_{j_k}| < \delta^{-1}$, so
	$$
	2\pi |\gamma_{r,k} - \omega_{r}\lambda_{j_k}| \leq \delta |1 - \exp(-2\pi i \omega_{r}\lambda_{j_k})|.
	$$
	
	Using the earlier inequality, we have
	\begin{align}\label{Ineq3ProofLemma3.4}
		|\gamma_{r,k} - \omega_{r}\lambda_{j_k}| \nonumber
		&< (2\pi)^{-1}\delta \max_{1 \leq r \leq m} \left|1 - \exp\left(-2\pi i\omega_{r}\lambda_{j_k} \right)\right| \nonumber \\
		&< (2\pi)^{-1}\delta \exp\left(-\epsilon j_k^{1/(2n\mu)}\right),
	\end{align}
	as $ k \to \infty $. 
	
	Now, we show that this estimate implies the failure of condition (i). Since
	$$
	|\gamma_{r,k}| - |\omega_r||\lambda_{j_k}| \leq |\gamma_{r,k} - \omega_r\lambda_{j_k}|,
	$$
	we deduce that
	$$
	|\gamma_{r,k}| \leq |\omega_r||\lambda_{j_k}| + |\gamma_{r,k} - \omega_r\lambda_{j_k}|.
	$$
	From \eqref{Ineq3ProofLemma3.4}, it follows that
	$$
	|\gamma_{r,k}| \leq C j_k^{M/(2n)},
	$$
	for some constant $ C > 0 $, and hence
	$$
	j_k \geq  C'\|\gamma_{k}\|^{2n/M}.
	$$
	
	Combining this with \eqref{Ineq3ProofLemma3.4}, we get for some $ \epsilon' > 0 $,
	$$
	|\gamma_{r,k} - \omega_r\lambda_{j_k}| \leq C_1 \exp\left\{-\epsilon'\left(\|\gamma_k\|^{1/(M\mu)} + j_k^{1/(2n\mu)}\right)\right\},
	$$
	for all $r\in \{1,\ldots,m\}$.
	
	Finally, since $ \sigma \geq M\mu $, we conclude that condition (i) fails, as desired.
\end{proof}

\textit{Proof of Proposition \ref{GHm=1}.}
From Theorem 3.11 in \cite{AviCap22} using Lemma  \ref{lt2}, we have that
$$
(a) \text{ or } (b) \Longrightarrow	GH. 
$$
Therefore, the sufficiency follows from \eqref{GHimpliesFGH}.

Now, if neither condition (a) nor (b) holds, it follows from \cite{AviCap22} that there exist singular solutions for $\mathcal{L}$. Specifically, we can obtain 
$$
u = \sum_{j \in \N} u_j(t) \varphi_j(x) \in \mathscr{U}_\mu(\mathbb{T} \times \R^n) \setminus \mathscr{F}_\mu(\mathbb{T} \times \R^n)
$$
such that $\mathcal{L} u \in \mathscr{F}_\mu(\mathbb{T} \times \R^n)$.

Moreover, these singular solutions are defined such that $\{u_j\} \subset \mathscr{G}^{\sigma}$ is uniformly bounded (see Theorems 3.6, 3.17, and Proposition 3.10). Therefore, by Lemma \ref{lemma_bound_seque}, we conclude that $u \in \mathscr{F}_\mu'(\mathbb{T} \times \R^n) \setminus \mathscr{F}_\mu(\mathbb{T} \times \R^n)$. Thus, $\mathcal{L}$ is not $\mathscr{F}_\mu$-globally hypoelliptic.
\qed

This result will be crucial in proving the following one.

\begin{proposition}\label{Lr-GH-implies-LL-GH}
	If $L_{r^*} = D_{t_{r^*}}+ c_{r^*}(t_{r^*})P(x,D_x)$ is  $\mathscr{F}_\mu$-globally hypoelliptic on $\TT^1 \times \R^n$, for some  ${r^*} \in \{1, \ldots, m\}$, then the system $\LL = (L_1, \ldots, L_m)$ is  $\mathscr{F}_\mu$-globally hypoelliptic on $\TT^m \times \R^n$.
\end{proposition}

\begin{proof}
	Let $F = (f_1, \ldots, f_m)$, where each coordinate $f_r \in \mathscr{F}_\mu(\TT^1 \times \R^n)$. By the definition of time-periodic Gelfand-Shilov spaces, for each $r \in \{1, \ldots, m\}$, there exists $\sigma>1$ such that $f_r \in \mathcal{S}_{\sigma, \mu}(\TT^1 \times \R^n)$ for every $r \in \{1, \ldots, m\}$. Possibly enlarging $\sigma,$ we can also assume that $c_r \in \mathcal{G}^\sigma(\TT^1)$ for $r=1,\ldots,m.$
	
	Now, assume that $u \in \mathscr{F}'_{\mu}$ is a solution of the system $L_r u = f_r \in \mathcal{S}_{\sigma, \mu}$ for $1 \leq r \leq m$, and consider the sequence of ordinary differential equations
	\begin{equation} \label{eq:ode}
		D_{t_r} u_j(t) + \lambda_j c_r(t_r) u_j(t) = f_{r,j}(t), \quad r = 1, \ldots, m,
	\end{equation}
	where $j \in \N$ and $t = (t_1, \ldots, t_m) \in \TT^m$.
	
	Since, by hypothesis, the operator $L_{r^*}$ is $\mathscr{F}_\mu$-globally hypoelliptic on $\TT_{r^*}^1 \times \R^n$, it follows from Proposition \ref{Z-infinity} that the set $\mathcal{Z}_{r^*}$ is finite. Therefore, for $r = r^*$ and for sufficiently large $j \in \N$, the equations in \eqref{eq:ode} have a unique solution, which can be expressed in the following two equivalent forms:
	\begin{equation}\label{solutions-1}
		u_j(t) = \dfrac{1}{1 - e^{-2\pi i \lambda_j c_{{r^*},0}}} 
		\int_{0}^{2\pi} \mathcal{H}_{{r^*},j}^-(t_{r^*}, \zeta)
		f_{r^*, j}(t_1, \ldots, t_{r^*} - \zeta, \ldots, t_m) \, d\zeta,
	\end{equation}
	and
	\begin{equation}\label{solutions-2}
		u_j(t) = \dfrac{1}{e^{2\pi i \lambda_j c_{{r^*},0}} - 1}
		\int_{0}^{2\pi} \mathcal{H}_{{r^*},j}^+(t_{r^*}, \zeta)
		f_{{r^*},j}(t_1, \ldots, t_{r^*} + \zeta, \ldots, t_m) \, d\zeta,
	\end{equation}
	where
	\begin{equation*}
		\mathcal{H}_{{r^*},j}^-(t_{r_0}, \zeta) =
		\exp \left(-i \lambda_j \int_{t_{r^*} - \zeta}^{t_{r^*}} c_{r^*}(\eta) \, d\eta\right),
	\end{equation*}
	and
	\begin{equation*}
		\mathcal{H}_{{r^*},j}^+(t_{r^*}, \zeta) = 
		\exp \left(i \lambda_j \int_{t_{r^*}}^{t_{r^*} + \zeta} c_{r^*}(\eta) \, d\eta\right).
	\end{equation*}
	We also denote
	$$
	\Theta_j^- = |1 - e^{-2\pi i \lambda_j c_{{r^*},0}}|^{-1}, \quad \text{or} \quad
	\Theta_j^+ = |e^{2\pi i \lambda_j c_{{r^*},0}} - 1|^{-1}.
	$$
	
	Let us assume that $\lambda_j>0$, for all $j$. For the remaining cases, see Remark \ref{what-to-do-if-lambda<0}.
	
	It follows from Proposition \ref{GHm=1} that either $b_{r^*} \not\equiv 0$ and $b_{r^*}$ does not change sign, or $b_{r^*} \equiv 0$ and $a_{r^*, 0}$ satisfies condition $\mathscr{D}_{\sigma,\mu}$ for all $\sigma \geq M\mu$. In both cases, we have 
	$$
	\sup_{t \in \TT^m, \ \zeta \in [0, 2\pi]} |\mathcal{H}_{r^*, j}^-(t, \zeta)| \leq C_r, \quad \text{or} \quad
	\sup_{t \in \TT^m, \ \zeta \in [0, 2\pi]} |\mathcal{H}_{r^*, j}^+(t, \zeta)| \leq C_r, \quad j \in \N.
	$$
	
	Moreover, for $b_{r^*, 0} \neq 0$, we have 
	\begin{equation}\label{b_0diff0}
		0 < A_1 \leq \Theta_j^- \leq 1, \quad \text{or} \quad 0 < B_1 \leq \Theta_j^+ \leq 1, \quad j \in \N.
	\end{equation}
	
	On the other hand, if $b_{r^*, 0} = 0$, for every $\epsilon > 0$, there exists $C_{\epsilon} > 0$ such that
	\begin{equation}\label{boundOmega}
		\Theta_j^- \leq C_{\epsilon} \exp \left(\epsilon j^{\frac{1}{2n\mu}}\right), \quad \text{as } j \to \infty,
	\end{equation}
	or equivalently, 
	\begin{equation}\label{boundTheta}
		\Theta_j^+ \leq C_{\epsilon} \exp \left(\epsilon j^{\frac{1}{2n\mu}}\right), \quad \text{as } j \to \infty.
	\end{equation}
	
	Next, let us analyze the behavior of $|\partial_t^{\alpha} u_j(t)|$ as $j$ tends to infinity. Assume $b_r \geq 0$, and consider expression \eqref{solutions-2}.
	
	By Theorem \ref{charac_partial_fourier-functions}, there exist constants $\epsilon_0 > 0$ and $C_1 > 0$ such that
	\begin{equation} \label{sbuff2}
		\sup_{t \in \TT^m} | \partial_t^{\gamma} f_{r,j}(t)| \leq C_1^{|\gamma| + 1} (\gamma !)^{\sigma} \exp \left(-\epsilon_0 j^{\frac{1}{2n\mu}} \right),
	\end{equation}
	for all multi-indices $\gamma$ and $j \in \N$.
	
	Using Leibniz formula and taking $\epsilon = \epsilon_0 / 2$, we obtain
	\begin{align*}
		|\partial_t^{\alpha} u_j(t)| & \leq \Theta_j^+ \sum_{\beta=0}^{\alpha} \binom{\alpha}{\beta} \int_0^{2\pi} |\partial_t^{\beta} \mathcal{H}_{r,j}^+(t_r)| \, |\partial_t^{\alpha -\beta} f_j(t - s)| ds \\
		& \leq C_2^{|\alpha| + 1} (\alpha !)^{\sigma} \, \Theta_j  \exp \left(-\frac{\epsilon_0}{2} j^{\frac{1}{2n\mu}} \right).
	\end{align*}
	
	Finally, if $b_{r,0} \neq 0$, we apply \eqref{b_0diff0}. On the other hand, if $b_{r,0} = 0$, we take $\epsilon = \epsilon_0 / 4$ in \eqref{boundTheta}. In any case, we obtain constants $\varepsilon$ and $C$ such that
	$$
	|\partial_t^{\alpha} u_j(t)| \leq C^{|\alpha| + 1} (\alpha !)^{\sigma} \exp \left(-\varepsilon j^{\frac{1}{2n\mu}} \right),
	$$
	which implies that $u \in \mathcal{S}_{\sigma,\mu}$.
\end{proof}

\begin{remark}\label{what-to-do-if-lambda<0}
	In the proof above, we assumed that all eigenvalues $\lambda_j$ are positive. In the general case, where some $\lambda_j$ may be negative, it suffices to partition the spectrum of $P$ into two subsets according to the sign of each $\lambda_j$. Then, apply the appropriate expression from \eqref{solutions-1} or \eqref{solutions-2} to ensure the decay of the solutions, paying attention to the signs that appear in the exponentials of these expressions.
\end{remark}

\subsection{Sufficiency of conditions (I) and (II)}

Assuming that condition $(I)$ of Theorem \ref{main_Theorem} holds, the $\mathscr{F}_\mu$-global hypoellipticity of $\LL$ follows directly from Propositions \ref{GHm=1} and \ref{Lr-GH-implies-LL-GH}. In other words, if for some $r$, with $1 \leq r \leq m$, the function $b_r$ is not identically zero and does not change sign, then $\LL$ is $\mathscr{F}_\mu$-globally hypoelliptic. Next, we will prove that condition $(II)$ of Theorem \ref{main_Theorem} implies $\mathscr{F}_\mu$-global hypoellipticity.

\begin{theorem}\label{cond-II-implies-GH}
	Assume that
	$$
	J = \left\{ r \in \{1, \dots, m\} : b_r(t_r) \equiv 0 \right\} = \{r_1 < \dots < r_\ell \} \neq \emptyset,
	$$
	and that $a_{J0} = (a_{r_1, 0}, \dots, a_{r_\ell, 0}) \in \R^{\ell}$ satisfies $\mathscr{D}_{\sigma,\mu}$ for all $\sigma\geq M\mu$. Under these conditions, the system $\LL$ is $\mathscr{F}_\mu$-globally hypoelliptic.
\end{theorem}	
\begin{proof}
	By Theorem \ref{The-Normal}, it suffices to prove that the corresponding normal form $\LL_0$ is $\mathscr{F}_\mu$-globally hypoelliptic. Moreover, by re-ordering the variables on the torus, we can assume that $J = \{1, \ldots, \ell\}$ and express the system $\LL_0$ as
	$$
	\LL_0 = 
	\left\{
	\begin{array}{rl}
		L_1 =& D_{t_1} + a_{1,0} P, \\
		\vdots \, & \\
		L_\ell =& D_{t_\ell} + a_{\ell,0} P, \\
		L_{\ell+1} =& D_{t_{\ell+1}} + (a_{{\ell+1},0} + ib_{\ell+1}(t_{\ell+1})) P, \\
		\vdots \, & \\
		L_m =& D_{t_m} + (a_{m,0} + ib_m(t_m)) P.
	\end{array}
	\right.
	$$
	
	Now, consider the reduced system
	$$
	\LL_0^{\ell} = 
	\left\{
	\begin{array}{l}
		L_1 = D_{t_1} + a_{1,0} P, \\
		\vdots \\
		L_\ell = D_{t_\ell} + a_{\ell,0} P.
	\end{array}
	\right.
	$$
	If $\ell=m$, the result follows directly from Proposition \ref{Prop_real_valued_coeff}. Therefore, we proceed by assuming $\ell < m$.

	Let us consider the decomposition $\TT^m = \TT^{\ell} \times \TT^{m-\ell}$ and write $(t,x) = (t',t'',x) \in \TT^m \times \R^n$, where $t' = (t_1, \ldots, t_{\ell}) \in \TT^{\ell}$ and $t'' = (t_{\ell+1}, \ldots, t_m) \in \TT^{m-\ell}$. 
	
	Given a solution $u \in \mathscr{F}'_{\mu}$ of the system
	$$
	L_r u = f_r \in \mathscr{F}_{\mu}, \quad r = 1, \ldots, m,
	$$
	we consider the Fourier series expansions of $u$ and $f_r$ with respect to the variables $(t',x)$, given by
	$$
	u(t,x) = \sum_{j \in \N} \sum_{\tau \in \Z^{\ell}} \widehat{u}_j(\tau,t'') e^{i t' \cdot \tau} \varphi_j(x),
	$$
	and
	$$
	f_r(t,x) = \sum_{j \in \N} \sum_{\tau \in \Z^{\ell}} \widehat{f}_{r,j}(\tau,t'') e^{i t' \cdot \tau} \varphi_j(x),
	$$
	for $1 \leq r \leq m$ and $(t,x) \in \TT^m \times \R^n$.
	
	By substituting these series into the equations $L_r u = f_r$ for $r = 1, \ldots, \ell$, we obtain
	$$
	(\tau_r + a_{r, 0} \lambda_j) \widehat{u}_j(\tau,t'') = \widehat{f}_{r,j}(\tau,t''),
	$$
	for each $\tau = (\tau_1, \ldots, \tau_\ell) \in \Z^\ell$, $j \in \N$, and $t'' = (t_{\ell+1}, \ldots, t_m) \in \TT^{m-\ell}$.
	
	Now, for each pair $(\tau, j) \in \Z^{\ell} \times \N$, we choose $s^* \in \{1, \ldots, \ell\}$ such that
	$$
	|\tau_{s^*} + a_{s^*, 0} \lambda_j| = \max_{1 \leq s \leq \ell} |\tau_s + a_{s, 0} \lambda_j|.
	$$

	Arguing as in the proof of Proposition \ref{Prop_real_valued_coeff} we have that the set of zeros of the symbol $\sigma_{\mathbb{L}^{\ell}_0}(\tau, j)$ for $(\tau, j) \in \Z^{\ell} \times \N$, is finite. 
	Then it follows that there exists a constant $K > 0$ such that $\sigma_{\mathbb{L}_{\ell}}(\tau, j) \neq 0$ for all $(\tau, j) \in \Z^{\ell} \times \N$ with $|(\tau, j)| \geq K$. Consequently, we can write
	\begin{equation}\label{est1-coeff-J-nonempty}
		\widehat{u}_j(t'', \tau) = (\tau_{s^*} + a_{s^*,0}\lambda_j)^{-1} \widehat{f}_{s^*, j}(t'', \tau),
	\end{equation}
	for all $(\tau, j) \in \Z^{\ell} \times \N$ such that $|(\tau, j)| \geq K$.
	
	In view of \eqref{seq-partial-coeff-funct}, there exist constants $\epsilon_0 > 0$ and $C > 0$ such that
	\begin{equation}\label{est2-coeff-J-nonempty}
		|\partial_{t''}^{\alpha} \widehat{f}_{s^*, j}(t'', \tau)| \leq C^{|\alpha|} (\alpha !)^{\sigma} \exp\left[- \epsilon_0 \left(\|\tau\|^{\frac{1}{\sigma}} + j^{\frac{1}{2n\mu}}\right)\right],
	\end{equation}
	for all $t'' \in \TT^{m-\ell}$, $\alpha \in \N_0^{m-\ell}$, and $(\tau, j) \in \Z^\ell \times \N$.
	
	Moreover, since $a_{J0} = (a_{r_1, 0}, \dots, a_{r_\ell, 0}) \in \R^{\ell}$ satisfies $\mathscr{D}_{\sigma,\mu}$ for all $\sigma \geq M\mu$, for the same $\epsilon_0 > 0$, there exists a constant $C_{\epsilon_0} > 0$ such that
	$$
	\max_{1 \leq s \leq \ell} |\tau_s - a_{s,0}\lambda_j| \geq C_{\epsilon_0} \exp\left[-\frac{\epsilon_0}{2} \left( |\tau|^{\frac{1}{\sigma}} + j^{\frac{1}{2n\mu}} \right)\right],
	$$
	for any $\tau = (\tau_1, \ldots, \tau_\ell) \in \Z^{\ell}$ and $j \in \N$. This implies
	\begin{equation}\label{est3-coeff-J-nonempty}
		|\tau_{s^*} + a_{s^*,0}\lambda_j|^{-1} \leq C_{\epsilon_0}^{-1} \exp\left[\frac{\epsilon_0}{2} \left(\|\tau\|^{\frac{1}{\sigma}} + j^{\frac{1}{2n\mu}}\right)\right],
	\end{equation}
	for any $(\tau,j) \in \Z^{\ell} \times \N$ such that $\sigma_{\mathbb{L}}(\tau,j) \neq 0$.
	
	Therefore, for all $t'' \in \TT^{m-\ell}$, $\alpha \in \N_0^{m-\ell}$, and $(\tau,j) \in \Z^\ell \times \N$, we have
	$$
	|\partial_{t''}^{\alpha}\widehat{u}_{j}(t'', \tau)| \leq C_{\epsilon_0}^{-1} C^{|\alpha|} (\alpha !)^{\sigma} \exp\left[-\frac{\epsilon_0}{2} \left(\|\tau\|^{\frac{1}{\sigma}} + j^{\frac{1}{2n\mu}}\right)\right].
	$$
	
	It follows from Theorem \ref{charac_partial_fourier-functions} that $u \in \mathscr{F}_{\mu}$.
\end{proof}

\section{Necessary conditions for hypoellipticity}\label{necessity-section}

Given the system $\LL = (L_1, L_2, \ldots, L_m)$, where
$$
L_r = D_{t_r} + c_r(t_r)P(x, D_x), \quad r = 1, \ldots, m,
$$
consider the associated time-independent system $\LL^0 = (L_1^0, L_2^0, \ldots, L_m^0)$, defined by
$$
L_r^0 = D_{t_r} + c_{r,0}P(x, D_x), \quad r = 1, \ldots, m,
$$
where $\displaystyle c_{r,0} = \frac{1}{2\pi} \int_{0}^{2\pi} c_r(s) \, ds$.

Before examining the necessity of conditions (I) and (II) in Theorem \ref{main_Theorem}, we give a preliminary result.

\begin{theorem}\label{L-GH-implies-Lo-GH}
	If the system $\LL$ is $\mathscr{F}_\mu$-globally hypoelliptic, then the system $\LL^0$ is also $\mathscr{F}_\mu$-globally hypoelliptic.
\end{theorem}
\begin{proof}
	Assume, by contradiction, that the system $\LL^0$ is not $\mathscr{F}_\mu$-globally hypoelliptic. By  Theorem \ref{GH-thm} there exists $\sigma$ for which the condition \eqref{Dio-cond-time-ind} fails or the set $\mathcal{N}_{\LL^0}$ is infinite. The latter possibility is excluded since in that case the set $\mathcal{Z}_{\LL}$ would be infinite, then since $\LL$ is $\mathscr{F}_\mu$-globally hypoelliptic, from Proposition \ref{Z-infinity} we would get a contradiction. If \eqref{Dio-cond-time-ind} fails for some $\sigma$, then according with Lemma \ref{lt2}, there is $\epsilon > 0$ and a sequence $(\tau_k, j_k) \in \Z^m \times \N$ such that 
	\begin{equation}\label{ineq_contra}
		0<	\max_{1 \leq r \leq m} \left|1 - \exp\left(-2\pi i \lambda_{j_k} c_{r,0}\right)\right| < k^{-1} \exp\left(-\epsilon \left(\|\tau_k\|^{1/\sigma} + j_k^{1/(2n\mu)}\right)\right),
	\end{equation}
	for all $k \in \N$.

	Note that, for each $k \in \N$, we can assume $\lambda_{j_k} > 0$ and, eventually excluding a finite number of $\lambda_{j_k}$, that $j_k \notin \mathcal{Z}_{\LL}$, where $\mathcal{Z}_{\LL}$ is defined by \eqref{Z_L-set}. As a matter of fact, since $\LL$ is $\mathscr{F}_\mu$-globally hypoelliptic, by Proposition  \ref{Z-infinity}, $\mathcal{Z}_{\LL}$ is finite.
	In particular, for each $k$, we may find $r \in \{1, \ldots, m\}$ such that $c_{r,0} \lambda_{j_k} \notin \Z$.
	
	Let $t^* = (t^*_1, \ldots, t^*_m) \in \TT^m$ be defined as
	$$
	\int_{0}^{t^*_r} b_r(\zeta) \, d\zeta \doteq \max_{t_r \in [0, 2\pi]} \int_{0}^{t_r} b_r(\zeta) \, d\zeta, \quad r = 1, \ldots, m.
	$$
	
	For each $r \in \{1, \ldots, m\}$, consider an interval
	$$
	I_r = [\alpha_r, \beta_r] \subset (0, 2\pi) \text{ such that } t^*_r \notin I_r,
	$$
	and let $\Phi_r \in \mathcal{G}^\sigma(I_r)$ be a real-valued cut-off function satisfying $0 \leq \Phi_r \leq 1$ and
	$$
	\int_{0}^{2\pi} \Phi_r(\zeta) \, d\zeta = 1.
	$$
	
	Now, if $c_{r,0} \lambda_{j_k} \in \Z$, we define
	\begin{equation}
		u_{r,k}(t_r) =  \exp\left[-\lambda_{j_k} \left(i \int_{0}^{t_r} c_r(\zeta) \, d\zeta + \int_{0}^{t^*_r} b_r(\zeta) \, d\zeta \right)\right] \label{u_rk-when-c_r0lambda_jk-is-integer}
	\end{equation}
	for $t_r \in [0, 2\pi]$.
	
	Observe that $|u_{r,k}(t^*_r)| = 1$ and $|u_{r,k}(t_r)| \leq 1$ for all $t_r \in [0, 2\pi]$. Moreover, we have
	$$
	D_{t_r} u_{r,k}(t_r) + \lambda_{j_k} c_r(t_r) u_{r,k}(t_r) = 0,
	$$
	for all $t_r \in [0, 2\pi]$.
	Moreover, applying Faà di Bruno's formula to obtain the following estimate:
	\begin{align*}
		|\partial_{t_r}^\ell u_{r,k}(t_r)| & \leq \sum_{\gamma=1}^\ell \sum_{\substack{\ell_1 + \ldots + \ell_\gamma = \ell \\ \ell_\nu \geq 1, \, \forall \nu}} |\lambda_{j_k}|^\gamma \frac{\ell!}{\ell_1! \cdots \ell_\gamma!} \left(\prod_{\nu=1}^\gamma \partial_t^{\ell_\nu-1}(c(t)-c(t-s)) \right) \\
		& \leq \sum_{\gamma=1}^\ell \sum_{\substack{\ell_1 + \ldots + \ell_\gamma = \ell \\ \ell_\nu \geq 1, \, \forall \nu}} |\lambda_{j_k}|^\gamma \frac{\ell!}{\ell_1! \cdots \ell_\gamma!} C^{\ell-\gamma+1}[(\ell -\gamma)!]^\sigma \\
		& \leq \sum_{\gamma=1}^\ell \sum_{\substack{\ell_1 + \ldots + \ell_\gamma = \ell \\ \ell_\nu \geq 1, \, \forall \nu}} C^{\gamma} j_k^{\gamma M/2n} \frac{\ell!}{\ell_1! \cdots \ell_\gamma!} C_1^{\ell-\gamma+1}[(\ell -\gamma)!]^\sigma \\
		&\leq C^{\ell +1} (\ell !)^{\sigma} j_k^{\ell M/2n},
	\end{align*}
	where we used the estimate $\lambda_{j_k} \leq C_1 j_k^{M/2n}$.
	
	On the other hand, if $c_{r,0} \lambda_{j_k} \notin \Z$, let $g_{r,k}$ be the periodic extension of the function
	$$
	t_r \in [0, 2\pi] \mapsto \left(1 - \exp(-2\pi i \lambda_{j_k} c_{r,0})\right) \exp\left(-i \lambda_{j_k} \int_{t^*_r}^{t_r} c_r(\zeta) \, d\zeta\right) \Phi_r(t_r).
	$$
	
	Thus, the unique $2\pi$-periodic solution of
	$$
	D_{t_r} u_{r,k}(t_r) + \lambda_{j_k} c_r(t_r) u_{r,k}(t_r) = g_{r,k}(t_r),
	$$
	is given by
	\begin{align*}
		u_{r,k}(t_r) &= \frac{i}{1 - \exp(-2\pi i \lambda_{j_k} c_{r,0})} \int_{0}^{2\pi} \exp\left(-i \lambda_{j_k} \int_{t_r-s}^{t_r} c_r(\zeta) \, d\zeta\right) g_{r,k}(t_r - s) \, ds\\
		&= i\int_{0}^{2\pi} \exp\left(-i \lambda_{j_k} \int_{t_r-s}^{t_r} c_r(\zeta) \, d\zeta\right) \exp\left(-i \lambda_{j_k} \int_{t^*_r}^{t_r-s} c_r(\zeta) \, d\zeta\right) \Phi_r(t_r - s) \, ds\\
		&=i \int_{0}^{2\pi} \exp\left(-i \lambda_{j_k} \int_{t^*_r}^{t_r} c_r(\zeta) \, d\zeta\right) \Phi_r(t_r - s) \, ds \\
		&=i \exp\left(-i \lambda_{j_k} \int_{t^*_r}^{t_r} c_r(\zeta) \, d\zeta\right) \int_{0}^{2\pi} \Phi_r(t_r - s) \, ds.
	\end{align*}
	In this case as well, we have 
	$$
	|u_{r,k}(t_r)| \leq \left|\int_{0}^{2\pi} \Phi_r(t_r - s) \, ds\right| = \left|\int_{t_r}^{t_r+2\pi} \Phi_r(\eta) \, d\eta\right| \leq 2\pi
	$$	
	and
	\begin{align*}
		|u_{r,k}(t^*_r)| &  = \left|\int_{0}^{2\pi} \Phi_r(t_r^* - s) \, ds\right|
		= \left|\int_{t^*_r}^{t^*_r + 2\pi} \Phi_r(\eta) \, d\eta\right| \\
		& = \int_{supp (\Phi_r)} \Phi_r(\eta) \, d\eta = 1,
	\end{align*}    
	whereas for the derivatives we get	
	\begin{align*}
		|\partial_{t_r}^\ell u_{r,k}(t_r)|& \leq\sum_{\alpha\leq \gamma} \binom{\gamma}{\alpha}\int_{0}^{2\pi} 				\left| \partial_{t_r}^{\alpha} \exp\left(-i\lambda_{j_k} 	\int_{t_r-s}^{t_r} c_r(\zeta)d\zeta\right)\right| \\
		& \times \left|\partial_{t_r}^{\gamma-\alpha} \left[\exp\left(-i\lambda_{j_k} \int_{t_k^{(r)}}^{t_r-s} c_r(\zeta)d\zeta \right) \Phi_r(t_r-s)\right]\right| ds \\
		&\leq C^{\ell +1}(\ell !)^{\sigma} j_k^{\ell M /2n}.
	\end{align*}
	
	If we now define, for each $k \in \N$,
	$
	u_k(t) \doteq \prod_{r=1}^{m} u_{r,k}(t_r)
	$
	and
	\begin{align*}
		f_{r,k}(t) &
		\doteq\left(D_{t_r} + c_r(t_r) \lambda_{j_k}\right) u_k(t) \\
		&=
		\begin{cases}
			\displaystyle \left(\prod_{s=1, s \neq r}^{m} u_{s,k}(t_s)\right) g_{r,k}(t_r), & \text{if } j_k \notin \mathcal{Z}_r, \\[4mm]
			0, & \text{if } j_k \in \mathcal{Z}_r,
		\end{cases}
	\end{align*}
	then clearly that $L_r u=f_r, r=1,\ldots, m.$ Moreover,
	$$
	|u_k(t^*)| = \prod_{r=1}^{m} |u_{r,k}(t_r^*)| = 1
	$$
	and
	\begin{align*}
		|u_k(t)| & = \prod_{r=1}^{m} \left| \exp\left(-i \lambda_{j_k} \int_{t^*_r}^{t_r} c_r(\zeta) \, d\zeta\right)\int_{0}^{2\pi} \Phi_r(t_r - s) \, ds\right|\\
		&\leq \left|\int_{0}^{2\pi} \Phi_r(t_r - s) \, ds\right|^m\\
		&\leq (2\pi)^m.
	\end{align*}
	for all $k \in \N$ and $t  \in \TT^m$. From the previous estimates, it follows that
	$$	u(t,x) =  \sum_{k \in \N} u_{k}(t) \varphi_{j_k}(x) \in \mathscr{F}_\mu'\setminus \mathscr{F}_\mu.$$
	To obtain a contradiction of the hypothesis, it is sufficient to show that
	$$	f_r(t,x) =  \sum_{k \in \N} f_{r, k}(t) \varphi_{j_k}(x) \in \mathscr{F}_\mu
	$$
	for every $r=1,\ldots,m.$

	Now, if $\gamma = (\gamma_1, \ldots, \gamma_m) \in \N_0^m$, we have
	$$
	|\partial_t^{\gamma} f_{r, k}(t)| \leq \left( \prod_{s=1, s \neq r}^{m} |\partial_{t_s}^{\gamma_s} u_{s,k}(t_s)| \right) |\partial_{t_r}^{\gamma_r} g_{r,k}(t_r)|,
	$$
	whenever $j_k \notin \mathcal{Z}_r$. In view of \eqref{ineq_contra}, and since
	$$
	\left|\exp\left(-i \lambda_{j_k} \int_{t_r^*}^{t_r} c_r(\zeta) \, d\zeta\right)\right| \leq 1,
	$$
	we can see that $|\partial_{t_r}^{\gamma_r} g_{r,k}(t_r)|$ satisfies an estimate of the form \eqref{deccoeff}.
	Therefore, $f \in \mathscr{F}_\mu$ and so the system $\LL$ is not $\mathscr{F}_\mu$-globally hypoelliptic.		
\end{proof}

In what follows, we consider the situation in which neither condition (I) nor (II) in Theorem \ref{main_Theorem} is satisfied and we construct singular solutions, which are ultradistributions $u \in \mathscr{F}_\mu' \setminus \mathscr{F}_\mu$ for which $L_r u \in \mathscr{F}_\mu$ for all $r = 1, \ldots, m$. In this way, we establish that $\LL$ is not $\mathscr{F}_\mu$-globally hypoelliptic. Notice that the fact that none of (I) and (II) are satisfied is equivalent to the occurrence of one of the following cases:
\\

\noindent
\textbf{Case 1:} The set $J$ in \eqref{J} is empty and every function $t_r \mapsto b_r(t_r)$, for $r =1,\ldots,m,$ changes sign.
\\

\noindent
\textbf{Case 2:} $J \neq \emptyset,$ say $J=\{r_1< \cdots< r_\ell \},$ and there exists $\sigma \geq M\mu$ such that the vector $a_{J0}= (a_{r_1,0}, \ldots, a_{r_\ell,0})$ does not satisfy $\mathscr{D}_{\sigma, \mu}$ and $b_r$ changes sign for every $r \notin J$. We devote the next two subsections to treat these two cases separately.

\subsection{Case 1}

The following lemma, whose proof can be found in \cite[Lemma 5.10]{AviGraKir18}, establishes crucial properties of periodic functions that change sign. These properties are instrumental for the construction of singular solutions presented in this subsection.

%

\begin{lemma}\label{lemma-change-sign}
	Let $b_r$ be a smooth, real-valued, $2\pi$-periodic function defined on $\R$. If $b_r \not\equiv 0$ then the following properties are equivalent:
	\begin{enumerate}
		\item[a)] $b_r$ changes sign;
		
		\item[b)] There exist $t_0^{(r)} \in \R$ and $t^{r,*}, t_{r,*} \in (t_0^{(r)}, t_0^{(r)} + 2\pi)$ such that
		\begin{align*}
			\mathcal{B}_{t^{r,*}}^{(r)}(\tau) & \leq 0, \quad \forall \tau \in (t_0^{(r)}, t_0^{(r)} + 2\pi], \\[2mm]
			\mathcal{B}_{t_{r,*}}^{(r)}(\tau) & \geq 0, \quad \forall \tau \in (t_0^{(r)}, t_0^{(r)} + 2\pi);
		\end{align*}
		
		\item[c)] There exist $t_0^{(r)} \in \R$ and partitions
		\begin{align*}
			& t_0^{(r)} < \alpha^*_{(r)} < \gamma^*_{(r)} < t^{r,*} < \delta^*_{(r)} < \beta^*_{(r)} < t_0^{(r)} + 2\pi,  \\[1mm]
			& t_0^{(r)} < \alpha_{*}^{(r)} < \gamma_{*}^{(r)} < t_{r,*} < \delta_{*}^{(r)} < \beta_{*}^{(r)} < t_0^{(r)} + 2\pi,
		\end{align*}
		and positive constants $c^{r,*}, c_{r,*}$ such that the following estimates hold:
		\begin{equation}
			\max \left\{ \mathcal{B}_{t^{r,*}}^{(r)} (\tau) \, ; \, \tau \in [\alpha^*_{(r)}, \gamma^*_{(r)}] \cup [\delta^*_{(r)}, \beta^*_{(r)}] \right\} < -c^{r,*}, \label{ch-sign-max1b}
		\end{equation}
		and
		\begin{equation}
			\min \left\{ \mathcal{B}_{t_{r,*}}^{(r)} (\tau) \, ; \, \tau \in [\alpha_{*}^{(r)}, \gamma_{*}^{(r)}] \cup [\delta_{*}^{(r)}, \beta_{*}^{(r)}] \right\} > c_{r,*}. \label{ch-sign-min1b}
		\end{equation}
	\end{enumerate}
\end{lemma}

We now proceed with Case 1. We assume that for each $r = 1, \ldots, m$, the function $b_r(t_r)$ is not identically zero and changes sign. Consequently, for each $r$, there exist points $t_r^+$ and $t_r^-$ in $[0, 2\pi]$ such that
$$
b_r(t_r^+) > 0 \quad \text{and} \quad b_r(t_r^-) < 0.
$$
For each $r = 1, \ldots, m$ and $\eta \in [0, 2\pi]$, we define the integral $\mathcal{B}^{\, r}_{\eta}: [0, 2\pi] \to \R$ by
$$
\mathcal{B}^{\, r}_{\eta}(\tau) = \int_{\eta}^{\tau} b_r(\zeta) \, d\zeta.
$$
The conditions outlined in Lemma \ref{lemma-change-sign} (particularly items b) and c)) allow us to select specific points $t^{r,*}, t_{r,*}$ and intervals which will be fundamental for our construction.

With the same notation as in Lemma \ref{lemma-change-sign}, consider the intervals
\begin{equation*}
	I^{r,*} \doteq [\alpha^{*}_{(r)}, \gamma^{*}_{(r)}] \cup [\delta^{*}_{(r)}, \beta^{*}_{(r)}] \quad \text{and} \quad
	I_{r,*} \doteq [\alpha_*^{(r)}, \gamma_*^{(r)}] \cup [\delta_*^{(r)}, \beta_*^{(r)}],
\end{equation*}
and choose functions $g^*_r, g_*^r, \psi^*_r, \psi_*^r \in \mathcal{G}^{\sigma}(\TT)$ such that
\begin{align*}
	& \text{supp}(\psi^*_r) \subset [0, 2\pi], \quad \psi^*_r|_{[\alpha^*_{(r)}, \beta^*_{(r)}]} \equiv 1, \\[2mm]
	& \text{supp}(g^*_r) \subset [\alpha^*_{(r)}, \beta^*_{(r)}], \quad g^*_r|_{[\gamma^*_{(r)}, \delta^*_{(r)}]} \equiv 1,
\end{align*}
and
\begin{align*}
	& \text{supp}(\psi_*^r) \subset [0, 2\pi], \quad \psi_*^r|_{[\alpha_*^{(r)}, \beta_*^{(r)}]} \equiv 1, \\[2mm]
	& \text{supp}(g_*^r) \subset [\alpha_*^{(r)}, \beta_*^{(r)}], \quad g_*^r|_{[\gamma_*^{(r)}, \delta_*^{(r)}]} \equiv 1.
\end{align*}

\begin{theorem}\label{change_sign}
	Suppose that for each $r \in \{1, \ldots, m\}$, the function $b_r(t_r)$ is not identically zero on any subinterval of $[0, 2\pi]$ and changes sign. Then the system $\mathbb{L}$ is not $\mathscr{F}_\mu$-globally hypoelliptic.
\end{theorem}

\begin{proof}
	Without loss of generality, assume that $\lambda_j > 0$, and consider the sequence of functions $\{u_{r,j}\}_{j\in\N}$ in $\mathcal{G}^{\sigma}(\TT)$ defined by
	\begin{equation*}
		u_{r,j}(t_r) = g^*_r(t_r) \exp\left[\lambda_j \psi^*_r(t_r) \left( \mathcal{B}_{t^{r,*}}^{(r)}(t_r) - i a_{r,0}(t_r - t^{r,*}) \right)\right].
	\end{equation*}
	In particular, if $t_r \in \text{supp}(g^*_r)$, we have
	\begin{equation*}
		u_{r,j}(t_r) = g^*_r(t_r) \exp\left[ \lambda_j \left( \mathcal{B}_{t^{r,*}}^{(r)}(t_r) - i a_{r,0}(t_r - t^{r,*}) \right)\right],
	\end{equation*}
	and since $\mathcal{B}_{t^{r,*}}^{(r)}(t_r) \leqslant 0$ on $I^{r,*}$, we obtain $e^{\lambda_j \mathcal{B}_{t^{r,*}}^{(r)}(t_r)} \leqslant 1$. Also, $|u_{r,j}(t^{r,*})| = 1$ for all $j \in \N$, and 
	\begin{equation*}
		\{u_{r,j}(t_r)\} \rightsquigarrow u_r \in \mathscr{F}'(\TT_r \times \R^n) \setminus \mathscr{F}(\TT_r \times \R^n).
	\end{equation*}
	
	For the sequence
	\begin{equation*}
		f_{r,j}(t_r) = -i {g^*_r}'(t_r) \exp\left[\lambda_j \psi^*_r(t_r) \left( \mathcal{B}_{t^{r,*}}^{(r)}(t_r) - i a_{r,0}(t_r - t^{r,*}) \right)\right],
	\end{equation*}
	we have $\text{supp}(f_{r,j}) \subset I^{r,*}$ for any $j \in \N$, implying that
	\begin{align*}
		\left| \partial_{t_r}^{\gamma_r} f_{r,j}(t_r) \right| 
		& \leq C^{\gamma_r+1} (\gamma_r!)^{\sigma+1} \exp\left(-\kappa_r j^{\frac{1}{2n\mu}}\right),
	\end{align*}
	for some $\kappa_r > 0$ and $\gamma_r \in \N_0$. In particular,
	\begin{equation*}
		D_{t_r} u_{r,j}(t_r) + (a_{r,0} + i b_r(t_r)) u_{r,j}(t_r) = f_{r,j}(t_r).
	\end{equation*}
	
	Now set $t^* = (t^{1,*}, \ldots, t^{m,*})$, and define
	\begin{equation*}
		u_j(t) = \prod_{r=1}^{m} u_{r,j}(t_r).
	\end{equation*}
	Since $|u_j(t^*)| = 1$ for all $j \in \N$, we obtain 
	$$u = \sum_{j \in \N} u_j \varphi_j \in \mathscr{F}_\mu' \setminus \mathscr{F}_\mu.$$
	
	Moreover, defining
	\begin{equation*}
		F_{r,j}(t) = D_{t_r} u_j(t) + (a_{r,0} + i b_r(t_r)) u_j(t) = f_{r,j}(t_r) \prod_{s=1, \, s \neq r}^{m} u_{s,j}(t_s),
	\end{equation*}
	we have, for all $\gamma = (\gamma_1, \ldots, \gamma_m) \in \N_0^m$,
	\begin{equation*}
		\partial_t^{\gamma} F_{r,j}(t) = \partial_{t_r}^{\gamma_r} f_{r,j}(t_r) \prod_{s=1, \, s \neq r}^{m} \partial_{t_s}^{\gamma_s} u_{s,j}(t_s).
	\end{equation*}
	
	Note that the derivatives $\partial_{t_s}^{\gamma_s} u_{s,j}(t_s)$ are bounded by powers of $\lambda_j$. Therefore, $F_r = \sum_{j \in \N} F_{r,j} \varphi_j \in \mathscr{F}_\mu$ and $L_r u = F_r$, implying that $\LL$ is not $\mathscr{F}_\mu$-globally hypoelliptic.
\end{proof}

\begin{remark}
	Theorem \ref{change_sign} extends to the following case: there exists an interval $[t_0^{k}, t_1^{k}] \subset [0, 2\pi]$ and $\delta_{k} > 0$, for some $k \in \{1, \ldots, m\}$, such that
	\begin{align*}
		b_k(t) & > 0, \ \forall t \in (t_0^{k} - \delta_k, t_0^{k}), \\
		b_k(t) & = 0, \ \forall t \in [t_0^{k}, t_1^{k}], \\
		b_k(t) & < 0, \ \forall t \in (t_1^{k}, t_1^{k} + \delta_k).
	\end{align*}
	Indeed, in this case, we may consider cutoff functions $g_0^{k}$ and $g_1^{k}$ such that
	\begin{align*}
		& \text{supp}(g_0^{k}) \subset [t_0^{k} - \epsilon, t_0^{k} + \epsilon] \ \text{ and } \ g_0^{k}|_{[t_0^{k} - \epsilon/2, t_0^{k} + \epsilon/2]} \equiv 1, \\[2mm]
		& \text{supp}(g_1^{k}) \subset [t_1^{k} - \epsilon, t_1^{k} + \epsilon] \ \text{ and } \ g_1^{k}|_{[t_1^{k} - \epsilon/2, t_1^{k} + \epsilon/2]} \equiv 1,
	\end{align*}
	for $\epsilon$ sufficiently small. Also, we define
	\begin{align*}
		B_{0,k}(t_k) & = \int_{t_0^k}^{t_k} b_k(\zeta) \, d\zeta, \ t_k \in \text{supp}(g_0^k), \\
		B_{1,k}(t_k) & = \int_{t_1^k}^{t_k} b_k(\zeta) \, d\zeta, \ t_k \in \text{supp}(g_1^k).
	\end{align*}
	Therefore, we may repeat the previous arguments by considering
	\begin{align*}
		u_{j,k}(t_k) &= g_0^k(t_k) \exp \left[ \lambda_j \left( B_{0,k}(t_k) - ia_{k,0} (t_k - t_0^k) \right) \right] \nonumber \\
		& \qquad + g_1^k(t_k) \exp \left[ \lambda_j \left( B_{1,k}(t_k) - ia_{k,0} (t_k - t_1^k) \right) \right].
	\end{align*}
\end{remark}

\begin{remark}
	Finally, if $\lambda_j < 0$, this construction can be adapted as previously described by defining the sequences in terms of the intervals $I_{r,*}$ and the corresponding cutoff functions $\psi_{*}^r$ and $g_{*}^r$.
\end{remark}

\subsection{Case 2}

In this subsection, we consider the case where
$$
J = \left\{ r \in \{1, \dots, m\} : b_r(t_r) \equiv 0 \right\} = \{r_1 < \dots < r_\ell \} \neq \emptyset,
$$
the vector $a_{J0} = (a_{r_1, 0}, \dots, a_{r_\ell, 0}) \in \mathbb{R}^{\ell}$ does not satisfy $\mathscr{D}_{\sigma,\mu}$ for some $\sigma\geq M\mu$ and $b_r$ changes sign for every $r \notin J$.	
By re-ordering the variables on the torus, we can assume that $J = \{1, \ldots, \ell\}$ and express the system as
$$
\LL = 
\left\{
\begin{array}{rl}
	L_1 =& D_{t_1} + a_{1,0} P, \\
	\vdots \, & \\
	L_\ell =& D_{t_\ell} + a_{\ell,0} P, \\
	L_{\ell+1} =& D_{t_{\ell+1}} + \left(a_{\ell+1,0} + ib_{\ell+1}(t_{\ell+1})\right) P, \\
	\vdots \, & \\
	L_m =& D_{t_m} + \left(a_{m,0} + ib_m(t_m)\right) P,
\end{array}
\right.
$$
where the functions $b_r$ change sign for $\ell + 1 \leq r \leq m$, and the vector $(a_{1,0}, \ldots, a_{\ell,0})$ satisfies the following condition:

There exist $\varepsilon > 0$ and  $(\tau_k, j_k) = (\tau_{k,1},\ldots,\tau_{k,\ell}, j_k)\in \Z^\ell \times \N$ such that
\begin{equation}\label{fincond1}
	0 < \|\widetilde{\sigma}_{\LL}(\tau_k, j_k)\| < \exp\left[-\varepsilon \left(\|\tau_k\|^{1/\sigma} + j_k^{1/(2n\mu)}\right) \right],
\end{equation}
\indent where \(\widetilde{\sigma}_{\LL}(\tau_k, j_k) = (\tau_{k,1} + a_{1,0}\lambda_{j_k}, \ldots, \tau_{k,\ell} + a_{\ell,0}\lambda_{j_k})\), and
$$
\|\widetilde{\sigma}_{\LL}(\tau_k, j_k)\| = \max_{1 \leq s \leq \ell} |\tau_{k,s} + a_{s,0}\lambda_{j_k}|.
$$

Recall the notation $t = (t', t'') \in \TT^m$, where $t' = (t_1, \ldots, t_\ell)$ and $t'' = (t_{\ell+1}, \ldots, t_m)$. Now define
$$
v_j(t') = 
\begin{cases}
	\exp(i t' \cdot \tau_k), & \text{if } j = j_{k}, \\
	0, & \text{otherwise},
\end{cases}
$$
and
$$
g_{r,j}(t') = 
\begin{cases}
	(\tau_{k,r}+ a_{r,0} \lambda_{j_{k}}) \exp(i t' \cdot \tau_k), & \text{if } j = j_{k}, \\
	0, & \text{otherwise},
\end{cases}
$$

Thus, we have
\begin{align*}
	v(t', x) &= \sum_{j \in \N} v_j(t') \varphi_j(x) \in \mathscr{F}_{\mu}'(\TT^{\ell} \times \R^n) \setminus \mathscr{F}_{\mu}(\TT^{\ell} \times \R^n), \\
	g_r(t', x) &= \sum_{j \in \N} g_{r,j}(t') \varphi_j(x) \in \mathscr{F}_{\mu}(\TT^{\ell} \times \R^n),
\end{align*}
and it follows that $L_r v = g_r$ for $r = 1, \ldots, \ell$.

On the other hand, by Proposition \ref{GHm=1}, each operator
$$
L_s = D_{t_s} + (a_{s,0} + ib_{s}(t_s)) P, \quad s = \ell+1, \ldots, m,
$$
is not $\mathscr{F}_\mu$-globally hypoelliptic on $\TT_{t_s} \times \R^n$ since $b_s$ changes sign for all $s = \ell+1, \ldots, m$. Moreover, from the proof of Theorem 3.17 in \cite{AviCap22}, we can construct a sequence $\omega_{j}(t_s)$ in the Gevrey class $\mathcal{G}^{\sigma}(\TT_{t_s})$ such that 
$$
\omega_s(t_s, x) = \sum_{j \in \N} \omega_{j}(t_s) \varphi_j(x) 
\in \mathscr{F}_{\mu}'(\TT_{t_s} \times \R^n) \setminus \mathscr{F}_{\mu}(\TT_{t_s} \times \R^n), \quad \ell + 1 \leq s \leq m,
$$
and $L_s \omega_s = h_s \in \mathscr{F}_\mu(\TT_{t_s} \times \R^n)$. Additionally, $\omega_{j}(t_s)$ is defined so that $|\omega_{j}(t_s)| \leq 1$, and there exists a point $\zeta_s \in \TT_{t_s}$ where $|\omega_{j}(\zeta_s)| = 1$ for all $j \in \N$.
Furthermore, all derivatives of $\omega_{j}(\zeta_s)$ are bounded by powers of $\lambda_{j}$.

Now, consider
$$
\Gamma_j(t'') = \prod_{s=\ell+1}^{m} \omega_{j}(t_s) \quad \text{and} \quad
u_j(t', t'') = v_j(t') \cdot \Gamma_j(t'').
$$
Then we claim that
$$
u(t', t'', x) = \sum_{j \in \N} u_j(t', t'') \varphi_j(x)
$$
is a singular solution of the system $\LL$. Indeed, by defining
$\zeta = (\zeta_{\ell+1}, \ldots, \zeta_m) \in \TT^{m-\ell}$, we have
$$
|u_{j_k}(0, \zeta)| = |v_{j_k}(0)| \cdot |\Gamma_{j_k}(\zeta)| = 1, \quad \forall k \in \N,
$$
implying that $u(t', t'', x) \in \mathscr{F}_{\mu}'(\TT^m \times \R^n) \setminus \mathscr{F}_{\mu}(\TT^m \times \R^n)$.

Now, for $r \in \{1, \ldots, \ell\}$, we define
$$
f_{r,j}(t', t'') = g_{r,j}(t') \cdot \Gamma_j(t''),
$$
and for $r \in \{\ell + 1, \ldots, m\}$,
$$
f_{r,j}(t', t'') = v_j(t') \cdot h_{r,j}(t_r) \cdot \prod_{s=1, s \neq r} \omega_j(t_s).
$$

It is straightforward to verify that $u(t', t'', x)$ satisfies
$L_r u = f_r$, where
$$
f_r(t', t'', x) = \sum_{j \in \N} f_{r,j}(t', t'') \varphi_j(x), \quad r = 1, \ldots, m.
$$

Next, we show that $f_r \in \mathscr{F}_{\mu}(\TT^m \times \R^n)$. Therefore, $u$ is a singular solution of the system $\LL$.

Let $\alpha = (\gamma, \beta) \in \N_0^{\ell} \times \N_0^{m-\ell}$. For $r \in \{1, \ldots, \ell\}$, we have
\begin{equation}\label{reg_f_1}
	\partial_t^{\alpha} f_{r,j}(t', t'') = \partial_{t'}^{\gamma} g_{r,j}(t') \cdot \partial_{t''}^{\beta} \Gamma_j(t''),
\end{equation}
and for $r \in \{\ell + 1, \ldots, m\}$,
\begin{equation}\label{reg_f_2}
	\partial_t^{\alpha} f_{r,j}(t', t'') = \partial_{t'}^{\gamma} v_j(t') \cdot \partial_{t_r}^{\beta_r} h_{r,j}(t_r) \cdot \prod_{s=1, s \neq r} \partial_{t_s}^{\beta_s} \omega_j(t_s),
\end{equation}
where $\beta = (\beta_{\ell+1}, \ldots, \beta_m)$.

From \eqref{fincond1}, there exist constants $\epsilon' > 0$ and $C_1 > 0$ such that
$$
\sup_{t' \in \TT^{\ell}} | \partial_{t'}^\gamma g_{r,j}(t') | \leq
C_1^{|\gamma|+1} (\gamma!)^{\sigma} \exp \left[-\epsilon' j^{\frac{1}{2n\mu}} \right].
$$

Moreover, there exists a constant $C_{\epsilon'} > 0$ such that
$$
\sup_{t'' \in \TT^{m-\ell}} | \partial_{t''}^\beta \Gamma_j(t'')| \leq C_{\epsilon'}^{|\beta|+1} (\beta!)^{\sigma} \exp \left[\frac{\epsilon'}{2} j^{\frac{1}{2n\mu}} \right].
$$

Hence, for $r \in \{1, \ldots, \ell\}$, it follows from \eqref{reg_f_1} that
$$
|\partial_t^{\alpha} f_{r,j}(t', t'')| \leq C^{|\alpha|+1} (\alpha!)^{\sigma} \exp \left[-\frac{\epsilon'}{2} j^{\frac{1}{2n\mu}} \right],
$$
implying that $f_r \in \mathscr{F}_{\mu}(\TT^m \times \R^n)$.
A similar argument applies to \eqref{reg_f_2}.

\section*{Acknowledgments}

The first and third authors thank the support provided by the National Council for Scientific and Technological Development - CNPq, Brazil (grants  305630/2022-9 and 316850/2021-7305630/2022-9, respectively). The second author has been partially supported by the Italian Ministry of the University and Research - MUR, 
within the PRIN 2022 Call (Project Code
2022HCLAZ8, CUP D53C24003370006).

\bibliographystyle{plain} 
\bibliography{references} 
\addcontentsline{toc}{section}{References}

\end{document}